\documentclass[oneside,reqno]{amsart}
\usepackage{amsfonts}
\usepackage{soul}

\usepackage{rotate}
\usepackage{amsmath,amssymb}
\usepackage[english]{babel}
\usepackage{graphicx}
\usepackage{color}

   \newcommand{\beq}{\begin{equation}}
   \newcommand{\eeq}{\end{equation}}
   \newcommand{\beqs}{\arraycolsep1.5pt\begin{eqnarray}}
   \newcommand{\eeqs}{\end{eqnarray}\arraycolsep5pt}
   \newcommand{\beqsn}{\arraycolsep1.5pt\begin{eqnarray*}}
   \newcommand{\eeqsn}{\end{eqnarray*}\arraycolsep5pt}

\newtheorem{thm}{Theorem}[section]
\newtheorem{rem}[thm]{Remark}
\newtheorem{cor}[thm]{Corollary}
\newtheorem{prop}[thm]{Proposition}
\newtheorem{lemma}[thm]{Lemma}
\newtheorem{defn}[thm]{Definition}
 \newtheorem{ex}[thm]{Example}
\numberwithin{equation}{section}

\newtheorem{pft}{}                              % 

\newcommand{\bspf}{\begin{pft}\begin{proof}{\bf Sketch of proof}.\ 
\rm}

\def\R{\mathbb {R}}
\def\N{{\mathbb N}}

\def\Q{{\Bbb Q}}

\newcommand{\ds}{\rightarrow}

\newcommand{\la}{\lambda}

\newcommand{\om}{\Omega}
\newcommand{\Om}{\Omega}
\newcommand{\f}{\varphi_\Omega}

\def\B{{\mathcal B}}

\def\L{{\mathcal L}}

\newcommand{\wi}{W^{1,\infty}(\Om)}

\def\supess{\mathop{\rm ess\: sup }}

\newcommand{\li}{{L^\infty}}

\newcommand{\conv}{{\rm co}}
\newcommand{\Lip}{{\rm Lip}\,}

\def\supess{\mathop{\rm ess\: sup }}

\def\to{\rightarrow}

\title{The role of intrinsic distances in the relaxation of $L^\infty$-functionals}
\author{Maria Stella Gelli}\address{M.S. Gelli, Dip. di Matematica, Universit\`a di Pisa,  Largo Pontecorvo 5, 
56127 Pisa (Italy)}
\author{Francesca Prinari}\address{F. Prinari, 
Dip. di Matematica e Informatica, Universit\`a di Ferrara, Via Machiavelli 35, 
44121 Ferrara (Italy)} 
 
\begin{document}
%\begin{abstract}{ We consider a  supremal functional of the form $$F(u)= \supess_{x \in \Omega} f(x,Du(x))$$ 
%where $\Omega\subseteq \R^N$ is a regular bounded open set,  $u\in \wi$ and  $f$ is  a Borel function.  Under a mild assumption on the sublevel sets  of $F$,   we give a description of the sublevel sets of its lower semicontinuous envelope  with respect to the weak$^*$ topology in terms of the level sets of suitable difference quotients.  This result is instrumental to show that the  lower semicontinuous envelopes of $F$  with respect to the  weak$^*$ topology, the weak$^*$ convergence and the uniform convergence  are level
%convex (i.e. they  have convex sub-level sets). The proof of these results relies both  on a deep  analysis  of  the intrinsic distances  associated to $F$ and on a careful use of variational tools such as $Gamma$-convergence. }
%%Moreover a comparison with the intrinsic distances associated to the lower semicontinuous envelopes of $F$ is given. 
%
%\end{abstract}. 
\begin{abstract}{ We consider a  supremal functional of the form $$F(u)= \supess_{x \in \Omega} f(x,Du(x))$$ 
where $\Omega\subseteq \R^N$ is a regular bounded open set,  $u\in \wi$ and  $f$ is  a Borel function. Assuming that the intrinsic distances $d^{\lambda}_F(x,y):= \sup \Big\{ u(x) - u(y): \, F(u)\leq \lambda \Big\}$
are locally equivalent to the euclidean one for every $\lambda>\inf_{\wi} F$, we give a description of the sublevel sets of the weak$^*$-lower semicontinuous envelope of $F$  in terms of the sub-level sets of the  difference quotient functionals $R_{d^\lambda_F}(u):=\sup_{x\not =y} \frac{u(x)-u(y)}{d^\lambda_F(x,y)}.
$
As a consequence  we prove that the relaxed functional of positive $1$-homogeneous supremal functionals coincides with $R_{d^1_F}$. Moreover, for a more general supremal functional $F$ (a priori non coercive), we prove that the sublevel sets of its relaxed functionals with respect to the  weak$^*$ topology, the weak$^*$ convergence and the uniform convergence are convex. The proof of these results relies both on a deep  analysis  of  the intrinsic distances  associated to $F$ and on a careful use of variational tools such as $\Gamma$-convergence.}
%Moreover a comparison with the intrinsic distances associated to the lower semicontinuous envelopes of $F$ is given. 

\end{abstract}

\maketitle

\bigskip\noindent{\it Mathematics Subject Classification (2000)}: 47J20, 58B20, 49J45.

\bigskip\noindent{{\it Keywords:} supremal functionals, level convex functions, intrinsic distances, relaxation, Calculus of Variations in $L^\infty$,  $\Gamma$- convergence.}

\bigskip
%%%%%%%%%%%%%%%%%%%%%%%%%%%%%%%%%%%%%%%%%%%%%%%%%%%%%%%%%%%%%%%%%%%%%%%%%%%%%%lower semicontinuous

%%%%%%%%%%%%%%%%%%%%%%%%%%%%%%%%%%%%%%%%%%%%%%%%%
%%%%%%%%%%%%%%%%%%%%%%%%%%%%%%%%%%%%%%%%%%%%%%%%%%

{\small \tableofcontents}

\section{Introduction}

In this paper we will consider functionals $F:  W^{1,\infty}(\Omega)\to \overline{\R}$ of the form
\beq\label{introf}
F(u)= \supess_{x \in \Omega} f(x,Du(x)), 
\eeq
where $f: \Omega \times \R^N \to  \overline{\R}$ is a Borel function and $\Omega\subseteq \R^N$ is a regular bounded open set. According to an established notation we will refer to energies in \eqref{introf} as  supremal or $L^\infty$-functionals on $W^{1,\infty}(\Omega)$ and we will use the term  supremand  to denote the function $f$   which represents the functional  \eqref{introf}.

\noindent A main topic in the study of supremal functionals %that has been carried out by different authors (mettere citazioni),
 is to give sufficient and necessary conditions on the supremand $f$ in order to have lower semicontinuity with respect to the  weak$^*$ topology  of $W^{1,\infty}(\Omega)$. Indeed, this topology is natural in order to study minimum problems since we have  compactness of minimizing sequences when $F$  satisfies a coercivity assumption.

\noindent A sufficient condition for  the sequential lower semicontinuity 
of a supremal functional with respect to the  weak$^*$ topology of $W^{1,\infty}(\om)$ 
has been shown by Barron,  Jensen and Wang in  \cite[Theorem 3.4]{BJW99}. It requires that for a.e. $x\in \Omega$ the function  $f(x,\cdot)$ is lower semicontinuous and  level convex, i.e.,  
$$
f(x,\theta \xi_1 + (1-\theta)\xi_2) \leq f(x,\xi_1) \vee f(x,\xi_2)
 \quad \forall \theta \in (0,1)\,  \ \forall\xi_1, \xi_2 \in \R^N.$$ 
%Equivalently   $f(x,\cdot)$ is  level convex  if and only if 
%for every $\lambda\in\R$ the sublevel set 
%$\big\{\xi\in\R^N\colon f(x,\xi)\le \lambda\big\}$ is convex. 
If $F$ is weak$^*$ lower semicontinuous,  in general the  vice-versa does not hold (\cite[Remark 3.1]{GPP}), due to the non-uniqueness of the supremand $f$ which represents $F$.
Note that if $f(x,\cdot)$ is level convex  then the supremal functional (\ref{introf}) is level convex, that is,  
for every $\la\in \R$   the sublevel set of $F$, denoted by 
\beq\label{Elaintro} E_\la=\Big\{ u\in \wi  : \, F(u) \le \lambda \Big\}, 
\eeq
is  convex. \\
\noindent In this paper we study the properties of the relaxed functional of  $F$ given by (\ref{introf}), i.e. the greatest lower semicontinuous functional $\Gamma_{\tau}(F)$ less than or equal to $F$ with respect to a fixed topology $\tau$  in $\wi$ chosen among 
the  weak$^*$, the sequential weak$^*$   and the uniform one (see Section \ref{prel}). This extensive analysis 
%emphasize the dependence on the topology 
relies on the fact that, in general, in absence of further coercivity hypotheses, $\Gamma_{\tau}(F)$ is affected by the choice of the topology, as  shown by Maggi and Gori  in \cite{MG}.  %Note that, under  suitable coercivity hypotheses on $F$, the relaxed functional $\Gamma_{\tau}(F)$ is independent of the topologies above. 

We note that, an explicit representation formula  of $\Gamma_{\tau}(F)$ in a supremal form has been established when $f$ is a continuous and coercive function (see \cite{P09}),  while the supremality of $\Gamma_{\tau}(F)$ is still unknown in the general case. Despite the lack of a representation result for $\Gamma_{\tau}(F)$,  in this paper we give a detailed description of its  sublevel sets by means of the intrinsic distance structures induced by $F$ (see Theorems \ref{relax-1} and \ref{relax0}).  As a result, for any $\tau$ quoted above,  we also show the level convexity of $\Gamma_{\tau}(F)$. %{\color{red} inserire anche il risultato di rappresentazione di funzionali 
%pos 1 omogenei non coercivi come funzionali quoziente, risultato che generalizza GPP per cui (forse) vale ancora %che il rilassato e' supremale se e solo se la distanza $d^1_F$ e' geodetica.}

%  the lower semicontinuous envelope of $F$ with respect to any topology $\tau$ quoted above.
%In particular,  if we denote by  $\tau$ any of these topologies, we get  that every  $\tau$-lower semicontinuous supremal functional is level convex. 

\noindent One of the main ingredients of our proofs is the introduction of the family of pseudo-distances $d_F^\lambda$ defined by 

\begin{equation}\label{dist}
d^{\la}_F(x,y):= \sup \Big\{ u(x) - u(y): \, F(u)\leq \la \Big\}
\end{equation}
for any $x,y\in \Omega$ and   $\la\in (\inf_{\wi} F, +\infty)$ (see  Section \ref{quotient}). 

These distances have been introduced by  De Cecco and Palmieri in the setting of Finsler metrics, in \cite{DCP1} and \cite{DCP2}, where they characterize the class of geodesic distances $d$ that satisfies the identity  $d=d^1_F$ for 
a suitable convex $1$-homogeneous supremal functional. 

 % i.e. 
%when $$f(x,Du(x))=\lim_{y\to x} \frac{|u(x)-u(y)|}{d(x,y)} \hbox{ a.e. }x\in \Omega.$$

 % Norris proved the small time asymptotics of the heat kernels to . were also studied by
 
 %Moreover, in \cite{DCP1} and \cite{DCP2} De Cecco and Palmieri considered the  metric derivative associated to a distance defined on a Lipschits manifold and locally  equivalent
%to the Euclidean one. They showed that in general the  metric derivative  cannot be expressed  almost everywhere as a square root of a quadratic form but belong to the wider class of the Finsler metrics.  %

In our  framework we consider a more general class of distance functions
which  are not symmetric or finite unless one requires additional hypotheses on $F$.  In  order to obtain a description of the sublevel sets of $\Gamma_{\tau}(F)$,  in Section \ref{quotient} 
we introduce the class of the  difference quotients $R_{d_F^\lambda}$ already considered   in \cite{GPP} and defined  by 
\begin{equation}\label{diffquot}R_{d_F^\lambda}(u)=\sup_{x,y\in \Omega, d_F^\lambda(x,y)\ne 0}\frac{u(x)-u(y)}{d_F^\lambda(x,y)}.
\end{equation}
Under the assumption that  the pseudo-distances $d^\la_F$ are metrically equivalent to the euclidean metric, that is, for every $\la>\inf_{\wi} F$,  $d^\la_F$ satisfies 
\begin{equation}
\label{metricequiv}
\alpha(\lambda)|x-y|\leq d^\la_F(x,y)\leq\beta(\lambda)|x-y|
\end{equation} 
with  $\alpha(\la),\beta(\la)$  strictly positive constants, 
in Theorem \ref{relax-1} we show that the 
sublevel sets of  $\Gamma_{\tau}(F)$ satisfy the identity
\beq\label{unione}
 \{ u\in\wi : \  \Gamma_{\tau}(F)(u)\leq\lambda  \}
=\bigcap_{\lambda'>\lambda} \{ u\in W^{1,\infty}(\Omega)\ : \ R_{d^{\lambda'}_F}(u)\leq 1   \}.
\eeq
As a consequence, under the assumption (\ref{metricequiv}), 
we get that  the relaxed functional $\Gamma_{\tau}(F)$ is independent of the topology  $\tau$ fixed above and it is level convex. 

Afterwards, exploiting Theorem  \ref{relax-1}, we provide  a representation  result for $\Gamma_{\tau}(F)$ when $F$ is  a positively $1$-homogeneous supremal functional satisfying \eqref{metricequiv}.  More in detail, in  Theorem \ref{generGPP} 
we represent $\Gamma_{\tau}(F)$ as the different quotient associated to the distance $d^1_F$, that is      
$$\Gamma_{\tau}(F)(u)=R_{d^1_F}(u) \quad \forall u\in\wi$$
 for any  $\tau$  topology among those quoted above. 
We underline that the class of difference quotients  strictly contains supremal functionals; under the assumptions of Theorem \ref{generGPP}, $\Gamma_\tau (F)$ can be represented in a supremal form if and only if the distance $d^1_F$ is \emph{intrinsic} according to the definition introduced by De Cecco and Palmieri \cite{DCP1,DCP2} (see Section 2 in \cite{GPP}).

Actually, the representation result in Theorem \ref{generGPP}  is a generalization of Theorem 3.5 in \cite{GPP}. Indeed, here we weaken the growth conditions 
%required in \cite{GPP} 
on the supremand $f$ representing $F$ and we remove its  continuity and symmetry hypotheses with respect to the gradient variable $\xi$.  

Note that  condition \eqref{metricequiv} has a key role in order to get the results quoted above. 
Hence, in order to clarify the setting  of validity of Theorems \ref{relax-1} and \ref{generGPP}, in  Theorem \ref{characterization} we characterize  the class of supremal functionals satisfying  \eqref{metricequiv}.  Roughly speaking, we prove that condition \eqref{metricequiv} holds  if and only if the sublevel sets of $F$ consist of functions with bounded gradients and $0\in \hbox{argmin }\Gamma_{\tau} (F) $  is a ''continuity'' point for $ \Gamma_{\tau} (F) $ along the affine functions $u_{\xi}:=\xi\cdot x$.

%The inequalities in  (\ref{metricequiv}) are related to the properties of the sub-level sets $E_{\la}$ given by (\ref{Elaintro}). More precisely the lower bound is verified, for instance, when $E_{\la}$ has non trivial interior part (e.g. if $f$ is a Carath\'eodory function) while 
%the upper bound  is equivalent to require that  $E_{\la}$ is bounded (i.e. $F$ is coercive).}
%{\color{red} Note that hypothesis  (\ref{dFbzero}) entails that for any $\lambda > \inf_{\wi} F$ the pseudo-distances $d^{\lambda}_F$ are metrically equivalent to the euclidean one.  
%Thus, in order to clarify the setting  of validity of Theorem \ref{relax-1}, in the following theorem  we provide a characterization of supremal functionals whose distances $d^{\lambda}_F$ satisfy (\ref{dFbzero}).}

Finally, in   Theorem \ref{relax0},  we study the relaxation problem for  a more general class of functionals $F$ with respect to those considered in Theorem \ref{relax-1}. In particular, we drop the boundedness assumption on the sublevel sets of $F$ thanks to an approximation result of $\Gamma_{\tau}(F)$  via $\Gamma$-convergence through a sequence of coercive ($\tau$-lower semicontinuous) functionals (see Proposition \ref{approxrelax}).  Moreover, %instead of requiring  a restrictive continuity assumption  on $f(x,\cdot)$, 
 we replace the condition that $0$ is a minimizer for $\Gamma_{\tau} (F)$ (together with  the continuity hypothesis above) with a mild hypothesis concerning the behaviour of $F$ along a suitable minimizing sequence. More precisely,  we assume  the  existence of a minimizing sequence for $F$ made up by  a sort of "upper semicontinuity" points, that is, $\exists (u_n)_{n\in \N}\subseteq \wi$ such that, set $u_\xi(x):=\xi\cdot x$, it holds   

\beq\label{accaomega0}\lim\limits_{n\to \infty} F(u_n)=\inf\limits_{\wi}F  
=\lim\limits_{n\to \infty}\limsup\limits_{\xi\to 0} F(u_n+u_\xi) 
\eeq
and we show that $\Gamma_{\tau}(F)$ is level convex when  $\tau$ is one of the topologies quoted above. 

Note that Theorem \ref{relax0}  applies to a wide class of supremal functionals $F$:  for instance,  assumption \eqref{accaomega0} is satisfied by  functionals $F$  whose supremand $f(x,\cdot)$  has a uniform modulus of continuity or  
satisfies the conditions
$$
 f(x,0)=0=\min\limits_{\xi\in \R^N}f(x,\xi) \hbox{ for a.e. }x\in \Omega \hbox{  and }
  \lim\limits_{\xi\to 0}  \supess\limits_{x\in \Omega}f(x,\xi)=0.
$$
In the particular case when $f$ does not depend on $x$, it is sufficient that  every sublevel sets of $f$ has non empty interior in order to get that $f$  satisfies \eqref{accaomega0} (for more details, see Remark \ref{hypo}).

As a byproduct of  Theorems \ref{relax-1} or  \ref{relax0},  we  obtain that the lower semicontinuous envelope of $F$   with respect to the weak$^*$ topology coincides with the sequential  weak$^*$ envelope. Note that, in general, the relaxed functional with respect to the  uniform convergence could be strictly  lower than the weak$^*$ envelope,  (see Example  2.2 in \cite{MG}).  
%for a  supremal functional $F$  not uniformly lower semicontinuous on $\wi$  represented by a  convex supremand $f$). 
Moreover, under the assumptions  of Theorem \ref{relax-1} or \ref{relax0},  we get also that if $F$ is $\tau$-lower semicontinuous  then $F$ is level convex. Note that with respect to the analogous result in  \cite{P09} (see Theorem 2.7), we drop the coercivity  and symmetry assumptions on $F$ and we do not require that $f$ is a Carath\'eodory function. 

We remark that a different notion of distance functions has been considered by Champion and De Pascale in \cite{CDP} 
in case the supremand $f(x,\xi)$ is level convex and coercive in $\xi$ and globally lower semicontinuous on $\Omega\times\R^N$ in order to establish comparison principles for absolute minimizers.  
In these hypotheses, using control theory formulas, they prove that their class of distances coincides with the family given by \eqref{dist}, main tool in our analysis. Moreover, the family of distances \eqref{dist} and the related different quotients  \eqref{diffquot} have been used by Davini and Ponsiglione in \cite{DP} in the setting of Finsler metrics, in order to study the closure  of two-phase gradient-constraints, in terms of 
$\Gamma$-convergence. 

We finally emphasize that the distance $d^1_F $ has been also considered in the particular case when $f(x,\xi)=\sqrt{\sum_{i,j} a_{ij}(x)\xi_i \xi_j}$ with $(a_{ij})$  an elliptic matrix in $\Omega$ or, more in general, when  $f$ is a Riemannian metric on a Lipschitz manifold $M$. In this case the  distance derived from the metric on $M$  has a relevant role  in the study of the heat flow associated to Dirichlet forms on $M$ (see \cite{N}), and, for a  smooth manifold,  $f$ coincides with the metric derivative  of the geodetic distance on $M$ (see \cite{DCP1, DCP2}). More recently, in  the  context of diffusion problem,  Koskela, Shanmugalingam and Zhou establish  when  the intrinsic differential and the local intrinsic distance structures coincide (see \cite{KSZ}).

\medskip

The paper is organized as follows:
Section \ref{prel} is devoted to some preliminary definitions and results concerning lower semicontinuous envelopes, level convex functionals and $\Gamma$-convergence; in Section \ref{mainres} we introduce the family of the pseudo-distances and their associated difference quotients and we state the main results of the paper. In Section \ref{quotient} we estabilish some key results about the family $\{d^\lambda_F\}_\lambda$ and its connection with the  family $\{d^\lambda_{\Gamma_\tau (F)}\}_\lambda$. These results will be instrumental for the proofs of the main results provided in Section \ref{proofs}. Eventually, Section \ref{repreresult} contains some additional results of interest. In particular we provide a representation result for $\Gamma_\tau (F)$ assuming its level convexity and we also address the problem of representing a weak$^*$ lower semicontinuous  functional $F$ by mean of a level convex supremand. 
Indeed, such a representation is crucial in problems involving supremal functionals as, for example, for  existence of absolute minimizers, in the homogenization problem, in the study of the $L^p$-approximation via $\Gamma$-convergence, in the characterization of the effective strength set in the context of electrical resistivity (see among others \cite{AP}, \cite{BJW}, \cite{BN}, \cite{CDPP02} and \cite{GPP}). 
Finally, in Section \ref{examples} we collect some interesting examples showing the optimality of some statements of the paper. %In Example 7.1 we provide two supremal functionals whose infimum is not supremal; 
In particular, in Example 7.2 we construct a supremal functional $F$ such that the  pseudo-distance $d^\lambda_{\Gamma_\tau (F)}$ associated to its relaxed functional $\Gamma_\tau (F)$ is different from  $d^\lambda_F$, for some values $\lambda >\inf_{\wi} F$; in Examples 7.3 and 7.5 we exhibit discontinuous supremands $f$ such that the associated supremal  functional $F$ still satisfies condition \eqref{accaomega}.  
%a supremal functional $F$ whose supremand $\f$ given by  \eqref{formulavarphi} does not represent the lower semicontinuous envelope of its localized version of $F$. 

%A functional $F:\wi\to \R\cup\{\pm \infty\}$ is {\bf level convex} if 
%$$
%F(\theta u_1 + (1-\theta)u_2) \leq F(u_1) \vee F(u_2)
%$$

%is lower semicontinuous functional with respect to the strong convergence in $L^\infty(\Omega)$ or with respect to the weak* convergence  in $\wi$, then  for every $\la\in \R$  the sublevel set $E_\la$ given  by (\ref{Ela}) is  convex.   In particular  $F$ is a level convex functional.

\section{Notations and preliminary results}\label{prel}

\noindent Throughout the paper we assume $\Omega$ to be a bounded open set in $\R^N$. We denote by ${\mathcal A}(\Omega)$  the family of all open subsets of $\Omega,$  and  by $\B_N$ the Borel $\sigma$-algebra of 
$\R^N$ (when $N=1$, we simply write $\B$). Moreover we write  $\|\cdot\|$ for 
the euclidean norm on $\R^N$,  $ B_r(x)$  for the open ball  $\{y\in\R^N \ : \|x-y\|<r \}$, and $\L^N$ for the Lebesgue measure in $\R^N$.  For any $a,b\in \bar \R$ we will denote  $a\wedge b:=\min \{a,b\}$ and $a\vee b:=\max \{a,b\}$. 
%For any couple of functions $f,g:\Omega\times \bar \R$ we denote by $f\wedge g$ and by  $f\vee g$ the functions  defined by $$(f\wedge g)(x)=f(x)\wedge g(x)  \quad \forall x\in \Omega$$ 
% $$(f\vee g)(x)=f(x)\vee g(x)   \quad \forall x\in \Omega.$$
If $\Omega $ is also connected, besides the euclidean one, it is possible to consider on  $\Omega$ the so called \emph{geodetic distance}, that is a distance containing the geometric features of both the open set and its boundary.   
More precisely, let  $\Gamma_{x,y}(\Om)$ be the set of Lipschitz  curves in $\Om$ with end-points $x$ and $y$. We define the \emph{geodetic distance } between two points $x,y\in \Omega$ as 
$$
|x-y|_\Om=\inf\{\mathcal{L}(\gamma) \ : \ \gamma\in \Gamma_{x,y} (\Om)\}\,
$$
where  $\mathcal{L}(\gamma)$ denotes  the length of the curve $\gamma$
with respect to the euclidean distance.
%$$
%\mathcal{L}(\gamma):= \sup \{ \sum_{i=1}^{k-1}  |\gamma(t_{i+1})-\gamma(t_{i})| \ : \  k\in\N,   \ 0=t_1<...<t_i<...<t_k=1 \}.
%$$

\noindent Note that  if  $\partial \Omega$ is Lipschitz then there exists a constant $C_{\Omega}>0$ such that 
\beq\label{normeequiv}|x-y|\le |x-y|_{\Omega}\leq  C_\Omega |x-y|.\eeq
We will use standard notations for Lebesgue and Sobolev spaces $L^p(\Omega), W^{1,p}(\Omega)$.  
We will also denote by $\Lip(\Omega)$ the space of the Lipschitz continuous functions on $\Omega$  and set 
$$\Lip_{\Omega}(u):=\sup_{x,y\in\Omega,\, x\not=y} \frac{|u(x)-u(y)|}{|x-y|}.$$ 
\noindent When $\Omega$ is bounded,   $ \Lip(\Omega) \subsetneq  W^{1,\infty}(\Omega)$ and  for any $u\in \Lip(\Omega)$ it holds $||Du||_{L^\infty (\Omega)}\leq \Lip_{\Omega}(u)$. 
In general, if  $\Omega$ is connected, then $||Du||_{L^\infty (\Omega)}$ coincides with the Lipschitz constant of $u$ with respect to the geodesic distance  as shown in the following Lemma. 
\begin{lemma}\label{Lipgeod}
Let $\Omega$ be a connected open set in $\R^N$. Then for any $u\in \wi$ it holds  
\beq\label{supnorma}
||Du||_{L^\infty (\Omega)}=\sup_{x,y\in\Omega, x\not=y} \frac{|u(x)-u(y)|}{|x-y|_\Omega}. 
\eeq
Moreover, if $\Omega$ is bounded and has Lipschitz continuous boundary then there exists a constant $c_\Omega>0$ such that 

\beq\label{normeLinfty}
||Du||_{L^\infty (\Omega)}\leq \Lip_{\Omega}(u)\leq c_\Omega ||Du||_{L^\infty(\Omega)}  \quad \forall\, u\in \wi .
\eeq
In particular,  $W^{1,\infty}(\Omega)\equiv \Lip(\Omega)$.  In addition, if  $\Omega$ is a convex set, then $||Du||_{L^\infty (\Omega)}=\Lip_{\Omega}(u)$. 
\end{lemma}
\proof It is well known that  $$|u(x)-u(y)|\leq ||Du||_{L^\infty (\Omega)} |x-y|_\Omega$$ for every $u\in \wi$ and for every $x,y\in \Omega$ (see, for example, \cite{Brezis}  Remark 7 in Chapter 9). 
The converse inequality in (\ref{supnorma}) can be established  by the fact that any $u\in\wi$ is differentiable almost everywhere in $\Omega$ and $Du$ coincides with the standard gradient of $u$.    Therefore, for every fixed $u\in \wi$ and for every ball $B\subset\subset \Omega$,  we have that 
$$
 \sup_{x,y\in\Omega,\, x\not=y} \frac{|u(x)-u(y)|}{|x-y|_\Omega}\geq \sup_{x,y\in B,\ x\not=y} \frac{|u(x)-u(y)|}{|x-y|_B} =\Lip _B (u)\geq ||Du||_{L^\infty (B)}.
$$
By passing to the supremum with respect to $B$ we obtain (\ref{supnorma}). Inequality \eqref{normeLinfty} follows by \eqref{normeequiv}. 
\qed

\bigskip

\noindent On $\wi$ we will consider different topologies:  the uniform topology  (denoted  by  $\tau_{\infty}$) induced by natural inclusion of $\wi\subseteq \li(\Omega)$  and  the weak* topology   ($w^*$ for shortly) inherited by $W^{1,\infty}(\Omega)$ as a (closed) subset of $L^{\infty}(\Omega)\times L^{\infty}(\Omega)$ endowed with the weak* topology. Moreover we will denote by $w^*_{seq}$  the  topology on $\wi$  induced by the $w^*$-convergence (see the next subsection).
 % as a dual space of $L^1(\Omega)\times L^1(\Omega)$. In particular, given $(u_n)_n, u\in W^{1,\infty}(\Omega)$ we write $u_n\wsto u$ in $\wi$ if $u_n\wsto u$ in  $L^\infty(\Omega)$ and $Du_n\wsto Du$ in $(L^\infty(\Omega))^N$. 

\subsection{Lower semicontinuous envelopes.} \label{topo}

Throughout the section  $(X,\tau)$ is a fixed topological space. For any set $B\subseteq X$ we denote by $\overline{B}^\tau$ its $\tau$-closure and we denote by  $\tau_{seq}$ the  topology on $X$  whose closed sets are the sequentially $\tau$-closed subsets of $X$. Note that $\tau_{seq}$ is in general strictly stronger than  $\tau$. 
%We provide the definition of lower semicontinuous envelope of a given function. 
\begin{defn}
Let $F:(X,\tau)\to \bar\R$ be a function. 

\noindent We say that $F$ is $\tau$-lower semicontinuous on $X$ (shortly $\tau$-l.s.c.) when for any $\lambda\in\R$ 
the sublevel set $\{x\in X\,|\, F(x)\le \lambda\}$ is $\tau$-closed. 
  
\noindent We say that $F$ is sequentially $\tau$-lower semicontinuous (shortly seq. $\tau$-l.s.c.) on $X$ if for any $x\in X$ and for any $(x_n)_n\subseteq X$ $\tau$-converging to $x$ we have 
$$
F(x)\le \liminf_{n\to +\infty} F(x_n).  
$$
\end{defn}

\begin{rem}  \label{property} {\rm 
 Note that
\begin{enumerate}
\item $F$  is sequentially $\tau$-lower semicontinuous if and only if $F$ is  $\tau_{seq}$-lower semicontinuous;

\item if $F:(X,\tau)\to \bar\R$  is $\tau\hbox{-lower semicontinuous}$ then $F$ is $\tau_{seq}$-lower semicontinuous;    

\item  the supremum of a family (also infinite) of $\tau\hbox{-lower semicontinuous}$  functions is still $\tau\hbox{-lower semicontinuous}$. Moreover, if $F,G$ are $\tau\hbox{-lower semicontinuous}$  then $F+G$ is $\tau\hbox{-lower semicontinuous}$.
\end{enumerate}}
\end{rem}

\begin{defn} Let $F:(X,\tau)\to \bar\R$. The {\sl lower semicontinuous envelope} (or {\sl relaxed function}) of $F$ is defined as 
\begin{equation}\label{envelope}
\Gamma_\tau(F) :=\sup\{G\,|\, G:(X,\tau)\to \bar\R \, ,   G \ \tau\hbox{-s.c.i.}\ \hbox{and } G\le F \hbox{ on } X\}. 
\end{equation}
\end{defn}
By Remark \ref{property}(2) it follows that $\Gamma_\tau (F)\le\Gamma_{\tau_{seq}} (F)$; if $(X,\tau)$ satisfies the first axiom of countability then  we have the following sequential characterization:
\def\totau{\stackrel{\tau}{\rightarrow}}
$$
\Gamma_\tau (F)(x)=\Gamma_{\tau_{seq}} (F)(x)=\min\big\{\liminf_{n\to +\infty} F(x_n)\,|\, x_n\totau x\big\}, 
$$
for any $x\in X$  (for a proof see \cite{Bu89} Proposition 1.3.3).  

\medskip

The following properties can be easily shown (for more details see \cite{DM93}, Chapters 3 and 6).
\begin{prop}\label{relaxsomma} Let $F:(X,\tau)\to \bar\R$, Then the following properties hold. 
\begin{enumerate}

\item [\rm (1)] $\Gamma_\tau(F)$ is $\tau$-lower semicontinuous;

\item [\rm (2)] $\inf_X F=\inf_X \Gamma_\tau(F)$;

\item [\rm (3)] for any $\tau$-continuous function $G:X\to\R$, we have $$\Gamma_\tau(F+G)=\Gamma_\tau(F)+G;$$ 

\item   [\rm (4)]  for any $c\in \R$, set $(F\vee c) (x):=F(x)\vee c$,    we have
\beq\label{fsupc} \Gamma_{\tau}(F\vee c)=\Gamma_{\tau}(F)\vee c;
\eeq

\item [\rm (5)]  if $X$ is a topological vector space and  $x_0\in X$, set  $G(\cdot):=F(\cdot+x_0)$, we have  
$$\Gamma_\tau(G)(\cdot)=\Gamma_\tau(F)(\cdot +x_0). 
$$

\end{enumerate} 
\end{prop}
In the following remark we clarify the relationship about the properties of lower semicontinuity of a functional $F:\wi\to \overline{\R}$ with respect to  the topologies
$w^*$, $w_{seq}^*$ and $\tau_{\infty}$.
\begin{rem} \label{lsc}{\rm Let $F:\wi\to \overline{\R}$ be a functional. 
Then the following properties hold. 
\begin{enumerate}
\item [\rm (1)]  If  $\forall \la\in \R$ there exists $r(\la)>0$ such that  \beq\label{weakcoerc1}E_{\la}=\{u\in \wi:\ F(u)\leq \la\}\subseteq \{ u\in \wi : \ ||Du||_{\li(\Omega)}\leq r(\la)\}\eeq then $$F  \  \hbox{ $w_{seq}^*$-l.s.c.}  \Longrightarrow \  F \ \tau_{\infty}\hbox{-l.s.c.};$$
indeed, thanks to \eqref{weakcoerc1}, any uniformly  convergent sequence $(u_n)_n$ in $\wi$ with
$\liminf\limits_{n\to\infty} F(u_n)<+\infty$ is also $w^*$-convergent in $\wi$;
\item [\rm (2)]  if  $\Omega$ has Lipschitz continuous boundary, by using standard immersion argument, it follows that $$F  \ \tau_{\infty}\hbox{-l.s.c.} \Longrightarrow F  \hbox{ $w_{seq}^*$-l.s.c.}.$$
Taking into account ($1$), if $\Omega$ has   Lipschitz continuous boundary  and $F$ satisfies \eqref{weakcoerc1} then
$$ F \hbox{  is $\tau_{\infty}$-l.s.c. } \Longleftrightarrow  F \hbox{ is  $w_{seq}^*$-l.s.c.} . $$
In particular  it holds
$$\Gamma_{\tau_{\infty}} (F)\equiv \Gamma_{w_{seq}^*}(F).$$ 

%\item  if $F(u,\Omega)\geq C\|Du\|_{L^{\infty}(\Omega)}$ for every $u\in \wi$ and is sequentially weakly* lower semicontinuous in $W^{1,\infty}( \Omega)$ then $F$ is also  lower semicontinuous with respect to the uniform convergence in $W^{1,\infty}( \Omega)$. 

\item [\rm (3)]  For what the $w^*$ topology is concerned, the following relation holds true without additional hypotheses on $F$ and $\Omega$  
$$ F \hbox{ is  $w^*$-l.s.c.}  \Longrightarrow  \ F \hbox{ is  $w_{seq}^*$-l.s.c.};$$
if  the sublevel sets $E_{\la}$ of $F$ are bounded in $\wi$, then the implication above can be reversed  
$$ F \hbox{ is  $w_{seq}^*$-l.s.c.}  \Longleftrightarrow  \ F \hbox{ is  $w^*$-l.s.c.}.$$ 
 Indeed on  bounded sets of $\wi$ the weak* topology induced by $L^\infty\times L^\infty$ is  metrizable. 
%In particular this result holds when $\Omega$ has   Lipschitz continuous boundary and  $F$ satisfies  the weaker assumption \eqref{weakcoerc}. 
\end{enumerate}
}
\end{rem}
\begin{rem}\label {MG1} {\rm  If we drop (\ref{weakcoerc1}), then the identity above could fails. Indeed in \cite[Example 2.2]{MG} Maggi and Gori  exhibit an example of a convex supremal functional $F$ that is  $w^*_{seq}$-l.s.c. (and  $w^*$-l.s.c.) but not $\tau_{\infty}$-l.s.c..   However, when $f(\cdot,\xi)$ is uniformly continuous, they show that the supremal functional $F$ represented by $f$ is $w^*$-l.s.c. if and only if 
 $F$ is $\tau_{\infty}$-l.s.c. (see \cite[Theorem 1.4]{MG}). }
 \end{rem}
 
 \subsection{$\Gamma$-convergence}
 
In order to introduce  the notion of $\Gamma$-convergence let $(X,\tau)$ be a  topological space and denote by $\mathcal{U}(x)$ the set of all open neighbourhoods of $x$ in $X$ .
\begin{defn} Let $F_n: X \ds \overline{\R}$ be a function for every $n\in \N$.
The {\sl $\Gamma(\tau)$-lower limit} and the {\sl $\Gamma(\tau)$-upper limit} of the sequence
  $(F_n)_{n\in \N}$  are the functions from $X$ into $ \overline{\R}$  defined by
$$ \Gamma(\tau)\hbox{-}\liminf_{n\to \infty} F_n(x) :=\sup_{U\in \mathcal{U}(x) \  }\liminf_{n\to \infty} \inf_{y\in U}F_n(y) $$
$$ \Gamma(\tau)\hbox{-}\limsup_{n\to \infty} F_n(x) :=\sup_{U\in \mathcal{U}(x) \  }\limsup_{n\to \infty} \inf_{y\in U}F_n(y) $$
If there exists a function $F:X\to \overline{\R}$  such that $\Gamma (\tau) \hbox{-}\liminf\limits_{n\to \infty} F_n= \Gamma(\tau)\hbox{-}\limsup\limits_{n\to \infty} F_n,$
 then we write $$F=\Gamma(\tau)\hbox{-}\lim_{n\to \infty} F_n$$ and we say that the sequence  $(F_n)_{n}$
$\Gamma (\tau)$-converges to $F$ or that $F$ is the $\Gamma(\tau)$-limit of $(F_n)_{n}$.

\end{defn}

In the following proposition we summarize some properties of the $\Gamma$-convergence useful in the sequel (see \cite{DM93} Proposition 6.8, Proposition 6.11, Proposition 5.7, Remark 5.5). 

\begin{prop}\label{gammaprop}  Let $F_n: X \ds \overline{\R}$ be a function for every $n\in \N$.
Then 
\begin{itemize}

\item[\rm (1)]  both the   $\Gamma(\tau)\hbox{-}\liminf\limits_{n\to \infty} F_n$ and  $\Gamma(\tau)\hbox{-}\limsup\limits_{n\to \infty} F_n$  are $\tau$-lower semicontinuous on $X$;
\item[\rm (2)]   the sequence  $(F_n)_n$
$\Gamma(\tau)$-converges to $F$ if and only if the sequence of the relaxed functions $(\Gamma_{\tau}(F_n))_{n\in \N}$
$\Gamma(\tau)$-converges to $F$; 

\item [\rm (3)] if $(F_n)_n$ is a not increasing sequence which  pointwise converges to $F$  then  $\Gamma(\tau)\hbox{-}\lim\limits_{n\to \infty} F_n= \Gamma_\tau (F)$. In particular if $F_n=F$ for every $n\in \N$ then  $\Gamma(\tau)\hbox{-}\lim\limits_{n\to \infty} F= \Gamma_\tau (F);$

\item [\rm (4)] if $(F_n)_{n\in \N}$ is an increasing sequence of $\tau$-lower semicontinuous functions which  pointwise converges to $F$  then  $\Gamma(\tau)\hbox{-}\lim\limits_{n\to \infty} F_n= F$.

\end{itemize}
\end{prop}

\subsection{Level convex functionals}

In the framework of supremal functionals,  level convexity plays the same main role as  convexity in the setting of  integral functionals.

\begin{defn}\label{levconvx} Let $(X,\tau)$ be a topological vector space. A function $F:X \to \overline{\R}$ is {\sl level convex} if 
$$
F(\theta x_1 + (1-\theta)x_2) \leq F(x_1) \vee F(x_2) \quad \forall \theta \in (0,1),\ \forall x_1, x_2 \in X$$ that is, for every $\la\in \R$  the sublevel set 
\beq
\label{Ela} E_\la=\Big\{ x\in X : \, F(x) \le \lambda \Big\}
\eeq 
is  convex.
\end{defn}
Note that the level convexity is stable under both  pointwise and $\Gamma$-convergence.
\begin{prop}\label{lcGammalimite} Let $(X,\tau)$ be a topological vector space and for every $n\in \N$  let $F_n:X\to \overline{\R}$ be  a  level convex function.  Then 
\begin{itemize}
\item [\rm (1)] the function $ F^{\#}(x)=\limsup \limits_{n\to\infty} F_n(x)$   is level convex;
\item  [\rm (2)] the function $F''(x)=\Gamma(\tau)$-$\limsup\limits_{n\to\infty} F_n$ is level convex. 
\end{itemize}

\end{prop}

\proof
$(1)$ Let $x_1,x_2\in X$ and $\theta \in (0,1)$.  Then $$ F^{\#}(\theta x_1 + (1-\theta)x_2) = \limsup_{n\to\infty}  F_n(\theta x_1 + (1-\theta)x_2) \leq \limsup_{n\to\infty}( F_n( x_1)\vee F_n(x_2)).$$
Since for every pair of real sequences $(a_n)_{n\in \N},(b_n)_{n\in \N}  $ we have that \beq\label{limsupinf}\limsup_{n\to\infty}( a_n\vee b_n)\leq (\limsup_{n\to\infty}a_n)\vee (\limsup_{n\to\infty} b_n),\eeq
 we get that
$$ F^{\#}(\theta x_1 + (1-\theta)x_2) \leq  (\limsup_{n\to\infty} F_n( x_1))\vee (\limsup_{n\to\infty} F_n(x_2))=F^{\#}(x_1)\vee F^{\#}(x_2).$$

$(2)$ Let $x_1,x_2\in X$, $\theta \in (0,1)$ and let $x:=\theta x_1 + (1-\theta)x_2$. Without loss of generality, assume that $F''(x_1)\vee F''(x_2)<+\infty$.
Since the map $(x,y)\mapsto \theta x + (1-\theta)y$
is continuous from $X\times X$ into $X$, then  for every $U\in \mathcal{U}(x)$  there exist $U_1\in \mathcal{U}(x_1)$ and  $U_2\in \mathcal{U}(x_2)$  such that $U$ contains the set $V:=\{ t y_1+(1-t)y_2\in X : \, y_1\in U_1,\ y_2\in U_2  \}.$
Then, by applying the level convexity of the functions $F_n$  we obtain that 
\begin{eqnarray*}\inf_{y\in U} F_n(y)&\leq& \inf_{y\in V}F_n(y) = \inf_{y_1\in U_1, y_2\in U_2} F_n( t y_1+(1-t)y_2)\\
 &\leq& \inf_{y_1\in U_1} \inf_{y_2\in U_2}( F_n( y_1)\vee F_n(y_2))=\inf_{y_1\in U_1}  F_n( y_1)\vee \inf_{y_2\in U_2}F_n(y_2)  .
\end{eqnarray*}
Hence, thanks to (\ref{limsupinf}), 
\begin{eqnarray*}\limsup_{n\to \infty} \inf_{y\in U}F_n(y) &\leq& \limsup_{n\to \infty} ( \inf_{y_1\in U_1}  F_n( y_1)\vee \inf_{y_2\in U_2}F_n(y_2) )\\
&\leq& (\limsup_{n\to \infty}  \inf_{y_1\in U_1}  F_n( y_1))\vee (\limsup_{n\to \infty}   \inf_{y_2\in U_2}F_n(y_2) ) .
\end{eqnarray*}
In particular 
$$\Gamma(\tau)\hbox{-}\limsup_{n\to \infty} F_n(x) \leq \Gamma(\tau)\hbox{-}\limsup_{n\to \infty} F_n(x_1)\vee \Gamma(\tau)\hbox{-}\limsup_{n\to \infty} F_n(x_2).$$
% Thanks to the sequential characterization of the $\Gamma$-convergence in a metric space, there exist $(x^i_n)_n \subset X$  such that $x^i_n\to x_i$ in $X$ for every $i\in \{1,2\}$ and  
%$$ F(x_i)=\lim_{n\to \infty} F_n(x^i_n).$$
%Then $$  F(\theta x_1+(1-\theta) x_2)\leq \liminf_{n\to \infty}F_n(\theta x_n^1+(1-\theta) x_n^2)$$
%$$\leq \liminf_{n\to \infty} F_n(x_n^1)\vee F_n(x_n^2) = \lim_{n\to \infty} F_n(x_n^1)\vee F_n(x_n^2)=F(x_1)\vee  F(x_2) .$$

\qed

\begin{rem} {\rm In general, given a sequence of level convex functions $(F_n)_{n\in \N} $, the functions $F'= \Gamma(\tau)\hbox{-}\liminf\limits_{n\to \infty}F_n$  and  $F_{\#}:=\liminf\limits_{n\to \infty}F_n$  are  not level convex. It is enough to consider the sequence $F_n(x) = (x - (-1)^n)^2$. In this case
$$F_{\#}(x)= F'(x)= (x+1)^2\wedge (x-1)^2.$$ }
%Esempio 7.3 in Dalmaso
\end{rem}

\bigskip

 Eventually we mention a further property that holds for level convex functions defined on $X'$ when $X$ is a separable Banach space. 

\begin{prop}\label{weakseq}  Let $X$ be a separable Banach space and let  $F:X'\to \overline{ \R}$. If the relaxed function $\Gamma_{w^*_{seq}} (F)$ is level convex (where $w^*$ stands for the weak* topology on $X'$), then 
$$
\Gamma_{w^*} (F)=\Gamma_{w^*_{seq}} (F).$$
In particular, if $F$ is  a level convex function, then 
$$F \hbox{ is $w^*$-lower semicontinuous} \Longleftrightarrow F  \hbox{ is  $w_{seq}^*$-lower semicontinuous}.$$
\end{prop}

\proof For any $\lambda\in \R$, the set  $E_\lambda = \{\Gamma_{w^*_{seq}} (F)\le \lambda \}$ is convex and sequentially weak* closed. Hence, by  applying  Banach-Dieudonne-Krein-Smulian Theorem (see \cite[Theorem 3.33]{Brezis}) we get that $E_\lambda$  is 
weak* closed. Therefore  $\Gamma_{w^*_{seq}} (F)$ is $w^*$-lower semicontinuous, which implies in turn that $\Gamma_{w^*_{seq}} (F)\leq \Gamma_{w^*} (F)$. The other inequality follows by Remark \ref{property}(2).\qed
\\

\begin{rem}{\rm  Note that, for a general   functional $F:\wi\to \bar\R$,  the lower semicontinuity with respect to one of the topologies $\tau_{\infty}, w^*, w^*_{seq}$ does not imply  the level convexity of $F$. Indeed it is enough to consider the characteristic function of the complement of any $\tau$-closed set $C\subseteq \wi$ that is not convex. For instance $C$ can be chosen as the union of two closed disjoint balls. }
\end{rem}

\section{Main results}\label{mainres}

%In this section $\Omega$ will be a connected bounded open subset of $\R^N$ with Lipschitz continuous boundary.   
Given a  supremal functional $F$  of the form 
\beq\label{sfzero}
F(u)= \supess_{x \in \Omega} f(x,Du(x)), 
\eeq 
where $\Omega$ is a bounded connected open set with Lipschitz continuous boundary, 
in  this section we  provide a description of the sublevel sets of   $\Gamma_{\tau} (F)$, where $\tau$ is one  of the topologies $\tau_{\infty}, w^*, w^*_{seq}$ defined in Section \ref{prel}, and we show that  $\Gamma_{\tau} (F)$ is a level convex functional.

We recall that by Lemma \ref{Lipgeod}   $||Du||_{L^\infty (\Omega)}$ coincides with the Lipschitz constant of $u\in \wi$ with respect to the  geodetic distance  $d_{\Omega}(x,y)=|x-y|_{\Omega}$. 
%:= \sup \Big\{ u(x) - u(y): \, ||Du||_{\infty}\leq 1 \Big\}
This means that,  when $f(x,\xi)=|\xi|$,   the  supremal functional (\ref{sfzero})
 can be represented as the  "difference quotient  functional"
\beq\label{quotientfunz}
F(u)=\sup_{x,y\in\Omega, x\not=y} \frac{u(x)-u(y)}{d_{\Omega}(x,y)}. 
\eeq
In addition, it can be also proved that the geodesic distance satisfies 
$$
%d_{\Omega}(x,y)=
|x-y|_{\Omega}=\sup \{u(x)-u(y): \, \|Du\|_{\infty}\le 1\}.
$$ 
Starting from this observation and exploiting the works \cite{DCP1, DCP2}, a representation formula of type  (\ref{quotientfunz})  has been proved in \cite{GPP} for the relaxation of $1$-homogeneous supremal functionals of type (\ref{sfzero}) with suitable growth conditions of linear type. 

%An analogous representation as different quotient has been proved in \cite{GPP} for $1$-homogeneous supremal functionals of type (\ref{sfzero}) with suitable growth conditions. 

Motivated by this result,  we  introduce a family $(d^{\lambda}_F)_{\la}$ of  pseudo-distances associated to the sublevel sets of a general  functional $F:\wi\to \bar \R$. With a slight abuse of notation we refer to this family as intrinsic distances associated to $F$. Hence we define the corresponding difference quotient  functionals and   we investigate their main properties. 

For any $\lambda>\inf_{\wi} F$ let $E_{\la}:=\{u\in\wi\, :\, F(u)\le \lambda\}$.   For any $(x,y)\in \Omega\times \Omega$ we set 
 \begin{equation}\label{dlambdaF}
d^{\la}_F(x,y):= \sup \Big\{ u(x) - u(y): \, u\in E_\la \Big\}. 
\end{equation}
Note that although $d^{\la}_F$ is in general not symmetric, by the definition of $d^{\la}_F$ it follows straightforward that 
$ \forall x,y,z\in \Omega
$ it holds 
$$
%\beq\label{triang} 
d^{\la}_F(x,y)\leq  d^{\la}_F(x,z)+d^{\la}_F(z,y). 
$$

\noindent We define the difference quotient functional (associated to $d^{\la}_F$) the functional $R_{d^\lambda_F} : \wi\to [0,+\infty]$ given
by
\beq\label{diffquotRd}
R_{d^\lambda_F}(u):=\sup_{x,y\in \Omega,\ d^\lambda_F(x,y)\ne 0}\frac{u(x)-u(y)}{d^\lambda_F(x,y)}.
\eeq  
For the sake of brevity, in the sequel we refer to $d^{\la}_F$ as 'distances' and to $R_{d^\lambda_F}$ as difference quotients.

We are now in position to state the main results of this paper.\\

\begin{thm}\label{relax-1} Let  $\Omega$ be a bounded connected open set with Lipschitz continuous boundary and let $F$ be  the   supremal functional \eqref{sfzero} represented by a  Borel function $f:\Omega\times \R^N\to \bar \R$.  
Assume that for any $\lambda > \inf_{\wi} F$   there exist positive coefficients $\alpha (\lambda), \beta(\lambda)>0$ such that for every $x, y\in \Omega$ 
\begin{equation}
\label{dFbzero}
\alpha(\lambda)|x-y|\leq d^\la_F(x,y)\leq\beta(\lambda)|x-y|_\Omega\,.
\end{equation} 
Then  $$\Gamma_{\tau_{\infty}} (F)\equiv \Gamma_{w^*}(F)\equiv \Gamma_{w_{seq}^*}(F)$$ 
and  the following identities 
%of the sublevel sets of  $\Gamma_{\tau} (F)$ 
hold for any $\lambda > \inf_{\wi} F$: 
\begin{eqnarray}\label{rapplevels}
 \{ u\in\wi : \  \Gamma_{\tau} (F)(u)\leq\lambda  \}&
=&\{u\in\wi : \  R_{d^{\lambda}_{\Gamma_\tau (F)}}(u)\leq 1  \}\\
\ &=& \bigcap_{\lambda'>\lambda} \{ u\in W^{1,\infty}(\Omega)\ : \ R_{d^{\lambda'}_F}(u)\leq 1   \}\notag
% \\
%\ &=&\color{red}{ \bigcap_{\lambda'>\lambda}\tau\hbox{-}\conv\Big(\{ u\in\wi : \  F(u)\leq\lambda  \}\Big)}\notag
\end{eqnarray}

where $\tau$ is any of  the topologies $\tau_{\infty}, w^*, w^*_{seq}$. In particular, $ \Gamma_{\tau} (F)$  is a level convex functional.   
%In addition  if $F$ is lower semicontinuous with respect to  one of the  topologies $\tau_{\infty}, w^*, w^*_{seq}$, then it also holds 
%\beq\label{rapplevels2}
% \{ u\in\wi : \  F(u)\leq\lambda  \}
%=\{ u\in W^{1,\infty}(\Omega)\ : \ R_{d^{\lambda}_F}(u)\leq 1   \}.   
%\eeq 

\end{thm}

 As a consequence of Theorem \ref{relax-1} we obtain the following representation result for the $\tau$-lower semicontinuous envelope  of a supremal functional depending only on the gradient and positively $1$-homogeneous. 
\begin{thm} \label{generGPP} Let  $\Omega$ be a bounded connected open set with Lipschitz continuous boundary and 
let $F$ be a   supremal functional represented by a  Borel function $f: \Omega\times\R^N\to \bar \R$. Assume that 
\begin{enumerate}
\item $F$  is positively $1$-homogeneous, that is, $F(\lambda u)=\lambda F(u)$  $\forall \lambda >0$
$\forall u\in \wi$;
\item $\exists \alpha, \beta>0$ such that $\forall x, y\in \Omega$ \ $\alpha|x-y|\leq d^1_F(x,y)\leq\beta|x-y|_\Omega$.
%the family $(d^{\lambda}_F)_{\la}$ satisfies (\ref{dFbzero}).
\end{enumerate}
Then    
\beq\label{reladiff}
\Gamma_{\tau}(F)\equiv R_{d^1_F}
\eeq
where $\tau$ is any of  the topologies $\tau_{\infty}, w^*, w^*_{seq}$.
\end{thm}

\begin{rem} \label{generGPP1}{\rm  Note that under the hypotheses of Theorem \ref{generGPP} it can be proved  that $\inf_{\wi} F=0$ and $\Gamma_\tau(F)(u)=\inf_{\wi} F=0$ if and only if $u=constant$. Moreover, Theorem  \ref{generGPP} provides a significant generalization of  the analogous result obtained in \cite{GPP}. Indeed in \cite{GPP} the  function $f(x,\xi):\Omega\times \R^N\to \bar \R$ representing the functional $F$ is assumed to be a Carath\'eodory function satisfying the conditions 
\begin{itemize}
\item[(i)] $f(x, \lambda \xi)=|\lambda| f(x, \xi)\quad \forall \la\in \R\,, \forall (x,\xi)\in \Omega\times\R^N$
\item[(ii)]$\exists C_1, C_2>0$ such that $ C_1|\xi|\le  f(x,  \xi)\le C_2 |\xi|\quad \forall (x,\xi)\in \Omega\times\R^N$.
\end{itemize}
Here we remove the continuity and symmetry hypotheses on  $f$  with respect to the variable $\xi$ and we weaken hypothesis (ii). More in detail, the linear growth condition on $f$ from below is replaced by the request that the sublevel sets  $\{(x,\xi)\,:\, f(x,\xi)\le \lambda\}$ are bounded in  $\Omega\times \R^N$, while the condition from above is fulfilled when $f$ is bounded on bounded sets of $\Omega\times \R^N$ (see Theorem \ref{characterization} for a deeper insight). 

Finally, we remark that, under the assumptions of Theorem \ref{generGPP}, it can be proved that $\Gamma_\tau (F)$ can be represented in a supremal form if and only if the distance $d^1_F$ is \emph{intrinsic} according to the definition introduced by De Cecco and Palmieri \cite{DCP1,DCP2} (see Proposition 2.5 in \cite{GPP}).}
%the representation result of $\Gamma_{\tau}(F)$ by a different quotient holds for functionals  }
\end{rem}

%{\color{red}  We emphasize that  Theorem \ref{relax-1}  above yields a  detailed description of the sublevel sets of $ \Gamma_{\tau} (F)$. However the metric equivalence of any $d^\lambda_F $ to the euclidean distance required by assumption \eqref{dFbzero}  is satisfied only by functionals $F$ whose sublevel sets  consist of functions with bounded gradients and such that $0\in \hbox{argmin }  \Gamma_{\tau} (F)$ (see Theorem \ref{characterization}). In  the following relaxation result we replace assumption \eqref{dFbzero}  with the \eqref{accaomega} and we cover a more general class of supremal functionals.}
%

We emphasize that,  at price  of hypothesis \eqref{dFbzero},  Theorem \ref{relax-1}  above yields a detailed description of the sublevel sets of $ \Gamma_{\tau} (F)$ that allows us, in Theorem \ref{generGPP}, to charactherize as a difference quotient the relaxed functional of $1$-homogeneous functionals. In the following theorem we replace assumption  \eqref{dFbzero}  with the a less restrictive hypothesis \eqref{accaomega} and prove  the level convexity of the relaxed functionals $\Gamma_{\tau}(F)$ for a more general class of supremal functionals.

%We emphasize that Theorem \ref{relax-1}  above yields a more detailed description of the sublevel sets of $ \Gamma_{\tau} (F)$ with respect to the next Theorem \ref{relax0}.  On the other hand hypothesis \eqref{dFbzero} in Theorem \ref{relax-1} is more restrictive  than hypothesis \eqref{accaomega} in Theorem \ref{relax0} that covers a more general class of supremal functionals. Indeed, the metric equivalence of any $d^\lambda_F $ to the euclidean distance, is satisfied only by functionals $F$ whose sublevel sets  consist of functions with bounded gradients and such that $0\in \hbox{argmin }  \Gamma_{\tau} (F)$ (see Theorem \ref{characterization}). 

\begin{thm}\label{relax0} Let  $\Omega$ be a bounded connected open set with Lipschitz continuous boundary and let $F$ be  the   supremal functional \eqref{sfzero} represented by a  Borel function $f:\Omega\times \R^N\to \bar \R$. Assume that 
%\noindent \eqref{accaomega} 
there exists $(u_n)_{n\in \N}\subseteq \wi$ such that, set $u_\xi(x):=\xi\cdot x$, it holds   
\beq\label{accaomega}\lim\limits_{n\to \infty} F(u_n)=\inf\limits_{\wi}F  
=\lim\limits_{n\to \infty}\limsup\limits_{\xi\to 0} F(u_n+u_\xi). 
\eeq
\noindent Then $ \Gamma_{\tau} (F)$  is a level convex functional when $\tau$ is one of  the topologies $\tau_{\infty}, w^*, w^*_{seq}$.  In particular $$\Gamma_{w^*}(F)\equiv \Gamma_{w_{seq}^*}(F).$$
\end{thm}

Finally, as a corollary of Theorems \ref{relax-1} and \ref{relax0}, we provide  the following result concerning the $\Gamma$-limits of sequence of supremal functionals.
\begin{cor} \label{stability} Let  $\Omega$ be a bounded connected open set with Lipschitz continuous boundary and 
let $(F_n)_n$ be a sequence of  supremal functionals 
\beq\label{sfzeron}
F_n(u)= \supess_{x \in \Omega} f_n(x,Du(x)), 
\eeq
represented by  Borel functions $f_n:\Omega\times \R^N\to \bar \R$. For every $ n\in \N$ assume that  $F_n$ satisfies either \eqref{dFbzero} or \eqref{accaomega}. If $F_n$ $\Gamma (\tau)$-converges  to a functional $F:\wi \to \bar \R$ when $\tau$ is one of  the topologies $\tau_{\infty}, w^*, w^*_{seq}$, then $F$  is level convex. 
\end{cor}

\begin{rem} \label{hypo}{\rm  
Note that hypothesis \eqref{accaomega} is satisfied in the following cases: 
\begin{itemize}
\item [-] when the  supremand  $f$ has  a uniform modulus of continuity in $\xi$, that is,  for any $M>0$ there exists 
  some modulus of continuity  $\omega_M$  such that 
  $$|f(x,\xi)-f(x,\eta)|\leq  \omega_M (|\xi -\eta|)\quad \hbox{for a.e. }x\in \Omega\,, \forall \xi,\eta \in B_M(0).
 $$ 
Such hypothesis has been already exploited in literature  to prove necessary conditions to the weak$^*$ lower semicontinuity of $F$ (see \cite{Pr08});
\item [-]  when  the supremand $f$ satisfies
$$
 f(x,0)=0=\min\limits_{\xi\in \R^N}f(x,\xi) \hbox{ for a.e. }x\in \Omega \hbox{  and }
  \lim\limits_{\xi\to 0}  \supess\limits_{x\in \Omega}f(x,\xi)=0, 
$$
(as in the model case when 
$0\leq  f(x, \xi) \le \alpha|\xi| \ \hbox{ for a.e. } x\in \Om\,, \forall \xi \in \R^N$ for some positive $\alpha$).   Indeed, it suffices to choose $u_n=0$ for $n\in \N$;
\item [-] when   $f=f(\xi)$  is such that  the  level set $H_{\la}(f):=\{ f(\xi)\leq \la\}$ has not empty interior for every $\la>\inf _{\R^N} f$. Indeed, let  $(\la_n)$ be such that $\la_n \to \inf_{\mathbb R^{N}}f$ and choose $\eta_n$  in the interior of 
 $H_{\la_n} (f).$ Then the sequence $u_n(x):=\eta_n\cdot x$ is a minimizing sequence and $\limsup_{\xi\to 0} F(u_n+u_\xi)=\limsup_{\xi\to 0} f(\eta_n+\xi) \leq \la_n$. Up to passing to the limit as $n\to +\infty$, we get  \eqref{accaomega}. 

\end{itemize}

}
\end{rem} 
%\begin{rem}\label{maggigori} {\rm We recall that, the relaxed envelopes $\Gamma_{\tau_{\infty}} (F)$   and  $\Gamma_{w_{seq}^*}(F)$ of a supremal functional $F$, may be different (see Remark \ref{MG1}). Nevertheless, Theorem \ref{relax0} ensures the level convexity of both these functionals.  We remark that we could not infer the level convexity of $ \Gamma_{\tau_{\infty}} (F)$  by the fact that $\Gamma_{w_{seq}^*}(F)$ enjoys this property since we  cannot establish  a priori that  $\Gamma_{\tau_{\infty}} (F)$ is a supremal functional.   
%}
%infatti se sapessimo che il rilassato uniforme e' supremale allora per la (2) del Remark 2.6 esso sarebbe anche semicontinuo seq debole e, applicando il Teorema 3.2, esso sarebbe level convesso perche' coinciderebbe con il suo rilassato debole seq. 
%\end{rem}
\begin{rem} \label{connessione}{\rm Note that Theorems~\ref{relax-1}  and \ref{relax0} could be generalized to more general bounded sets $\Omega$ in $\R^d$ by requiring that any connected component $\Omega_i$ of $\Omega$ is Lipschitz regular  with the constants $ c_{\Omega_i}$ appearing in the  inequality \eqref{normeequiv} and involving the geodetic distances associated to any connected component $\Omega_i$ are equibounded  (and  for Theorem \ref{relax0} a local version of  hypothesis \eqref{accaomega} must be satisfied). Indeed, the proof follows by arguing separately on each connected component. 
}
\end{rem}
\bigskip

%We refer to Proposition \ref{altoprop} for a sort of characterization of functionals satisfying \eqref{dFbzero}.

%: indeed, the right-hand side inequality is equivalent to the boundedness of the sublevel sets $E_\lambda$ of $F$, while the left-hand side inequality holds true provided each $E_\lambda$ contains a dense sets of 'directions', this last condition being accomplished when $F$ is strongly continuous in its minimum point along affine perturbations.  

\section{Intrinsic distances and related results.}\label{quotient}

In this section  we prove some key results instrumental for the proofs of Theorems \ref{relax-1} and \ref{relax0}. According to the hypothesis  (\ref{dFbzero}) of Theorems \ref{relax-1}  in the sequel  we assume  that   the  pseudo-distances $ d^\la_F$ are metrically equivalent to the euclidean distance on $\Omega$. More in detail, for every  $\lambda >\inf_{\wi} F$   $ d^\la_F$ satisfies 
$$
\alpha(\lambda)|x-y|\leq d^\la_F(x,y)\leq\beta(\lambda)|x-y|_\Omega\, \qquad \forall x,y\in \Omega
$$
for positive $\alpha(\lambda), \beta(\lambda)$.  
Note that the the validity of the right hand-side of the  inequality above ensures that  the (non symmetric)  distance $ d^\la_F$ is non-degenerate, that is  
$$d^\la_F(x,y)=0 \iff x=y.$$
In particular, for every $u\in \wi$  and for every  $\lambda >\inf_{\wi} F$, we can rewrite the different quotient as 
 $$
R_{d^\lambda_F}(u)=\sup_{x,y\in \Omega,\; x \neq y}\frac{u(x)-u(y)}{d^\lambda_F(x,y)} 
$$
 and, by using Lemma \ref{Lipgeod}, we also get that 
\beq\label{stimagrad}
 \frac 1 {\beta (\la)} ||Du||_{\infty}\leq R_{d^\lambda_F}(u)\leq \frac 1  {\alpha(\la)C_\Omega } {||Du||_{\infty}}
\eeq
where $C_\Omega$ is the constant in \eqref{normeequiv}. 
First of all we state some comparison results inherited straightforward by the definition that will be useful in the sequel. \\

\begin{itemize}
\item[(i)] 
 if $F, G:\wi\to \bar\R$ and $F\le G$, then $\forall \lambda >\inf_{\wi} G$ it holds 
\beq\label{dis1}
d^{\la}_G(x,y)\le d^{\la}_F(x,y) \quad \forall x,y\in\Omega,
\eeq
and 
\beq\label{dis2}
R_{d^{\la}_F}(u)\le R_{d^{\la}_G}(u) \quad \forall u\in\wi;
\eeq

\item[(ii)]   if $\lambda\le \lambda '$, then $\forall \lambda >\inf_{\wi} F$ it holds
\beqsn
d^{\la}_F(x,y)\le d^{\lambda'}_F(x,y) \quad \forall x,y\in\Omega, \\
R_{d^{\lambda'}_F}(u)\le R_{d^{\la}_F}(u) \quad \forall u\in\wi.
\eeqsn
\end{itemize}

\medskip

The following proposition follows the same guidelines of Lemma 3.4  in \cite{GPP} up to some suitable changes necessary to deal with the lack of symmetry and linear growth of $F$.  

\begin{prop}\label{keylemma} Let  $\Omega$ be a bounded connected open set with Lipschitz continuous boundary and let $F:\wi\to \bar \R$ be a  supremal functional of the form $F(u)= \supess_{x \in \Omega} f(x,Du(x))$
 represented by  a Borel function  $f:\Om \times \R^N \to \bar \R$. 
Assume that  $ (d^\la_F)_{\la}$ satisfies (\ref{dFbzero}). %for every  $\lambda >\inf_{\wi} F$.
 Then  for  any $v\in \wi$ such that $R_{d^\la_{F}}(v) < 1$, there exists
a sequence $\{v_n\}\subset\wi$ converging to $v$ in $\li$ with $F(v_n) \le \la$. 
\end{prop}

\proof Note that for any $x,y\in \Omega$   it holds 
\begin{equation}\label{emme} 
d^\la_F (x,y)=d^{\la-m}_{F-m}(x,y)  \quad \forall m\in \R, \quad\forall \lambda >\inf_{\wi} F.
\end{equation}  
We claim that, without loss of generality, we may confine to prove  the thesis for $\lambda$ positive. Indeed, let $\lambda > \inf_{\wi} F$ and $v\in \wi$ such that  $R_{d^\la_{F}}(v) < 1$ be fixed  and  select $\lambda'\in\R$ such that $\inf_{\wi} F<\lambda'<\lambda$. Then, thanks to \eqref{emme} with $m=\lambda'$  we have 
 $R_{d^{\lambda-\lambda'}_{F-\lambda'}}(v)=R_{d^\la_{F}}(v) < 1$. If we find $\{v_n\}\subset\wi$ converging to $v$ in $\li$ with $(F-\lambda')(v_n) \le \lambda-\lambda' $ then we also deduce that  $F(v_n) \le \la$. Since $\lambda -\lambda'>0$, we may then assume directly $\lambda$  to be strictly positive. Up to replace $F$ with $F/\lambda$  it is not restrictive to treat the case  $\lambda=1$ and drop the dependence on $\lambda$ in the notation of the associate distance and of $\alpha(\lambda), \beta(\lambda)$. 

\noindent Let  $w\in W^{1,\infty}(\Omega) $ satisfying $F(w)\leq 1$   then, by definition of  $ d_F(x,y)$, for $x,y\in \Omega$ we have  that 
$$
w(x)-w(y)\le d_F(x,y)\le \beta|x-y|_\Omega
$$  
and, switching the role of $x, y$, we deduce  
$$
|w(x)-w(y)|\le \beta|x-y|_\Omega.
$$
Thanks to Lemma \ref{Lipgeod} this implies \begin{equation}\label{gradgeod}
\|D w \|_{\infty} \leq  \beta.
\end{equation}

 Now let us fix $v\in \wi$ such that $R_{d_{F}}(v) < 1$
%and the regularity of $\partial \Omega$ implies  (\ref{lip}).
and a positive radius $r>0$. Let $0<\theta<1$ such that 
$R_{d_F}(v)=1-\theta < 1$. Then,  thanks to the assumption (\ref{dFbzero}),
for every $x$, $y\in \Omega$ with $|x-y|=r$
 we have that
 $$v(y)-v(x)\leq (1-\theta)d_F(y,x)=d_F(y,x)- \theta d_F(y,x)\leq d_F(y,x)-\alpha\theta|x-y|<d_F(y,x)-r\theta\alpha .$$
 In particular 
\begin{equation}\label{stima1} v(y)-v(x) < d_F(y,x)-r\theta\alpha.
 \end{equation}
%where  $\gamma=r\theta\alpha $.

Let us fix $0<\varepsilon<r\min\{\frac{1}{3}\theta\alpha, \beta+\|Dv\|_{\infty}\}$. For every $x\in \Omega$ 
and for every $y\in \partial B_r(x) \cap \Omega$, by the definition of $d_F$ there exists
a function $w_r^{x,y}\in \wi$ such that
\vskip4pt
\begin{itemize}
\item[1)]
$F(w_r^{x,y},\Omega) \le 1;$
\vskip4pt
\item[2)]
$ w_r^{x,y}(y)\ge w_r^{x,y}(x) + d_F(y,x) - \varepsilon;$
\vskip4pt
\item[3)]
$ w_r^{x,y}(x) = v(x)$;
\end{itemize}
\vskip4pt
the third property being possible thanks to the translation invariance of the first two.

By properties $2)$, $3)$ and by (\ref{stima1}), for every $y\in \partial B_r(x)\cap \Omega$ 
we have that 

$$
w_r^{x,y}(y)\ge w_r^{x,y}(x) + d_F(y,x) - \varepsilon= v(x) + d_F(y,x) -\varepsilon> v(y) + r\theta\alpha
  -\varepsilon
$$
that is 
$$
w_r^{x,y}(y)- v(y) > r\theta\alpha  -\varepsilon
$$
for every $y\in \partial B_r(x)\cap \Omega$.
Thanks to  (\ref{gradgeod}),  
%the functions $w_r^{x,y}- v$ are $L$-Lipschitz continuous.
for  $\delta=\frac {\epsilon }{\beta+\|Dv\|_{\infty}}r$,  we have that  for every $z\in \Omega$ such that $|z-y|_\Omega\le\delta$ we have that 
$$w_r^{x,y}(z)- v(z)> w_r^{x,y}(y)- v(y)- (\beta+\|Dv\|_{\infty}) |z-y|_\Omega\geq r\theta \alpha  -2\varepsilon,   $$
that is 
\begin{equation}\label{stima*}
w_r^{x,y}(z)> v(z)+\varepsilon \qquad \forall y \in \partial B_r(x)\cap\Om, \forall z\in \Omega :\, |z-y|_\Omega\le
\delta .
\end{equation}
Moreover, since $ w_r^{x,y}(x) = v(x),$ we have that for every $z\in \Omega$ such that $|z-x|_\Omega<\delta $
\begin{equation}\label{stima**}w_r^{x,y}(z) \leq  v(z)+ (\beta+\|Dv\|_{\infty}) |z-x|_\Omega< v(z)+\epsilon.  
\end{equation}
Note that, since $\delta \le \frac 13 r$, we have that $\forall x,y\in \Omega$ such that $|x-y|=r$ it holds $B_\delta(x)\cap B_\delta(y)\cap \Omega=\emptyset$. 
Moreover  the family $$\{ B_{\delta}(y): \ y\in \Omega\cap \partial B_r(x)\}$$ is an open covering of  the pre-compact set $\partial B_r(x)\cap  \Omega$. 
For every $x\in \Omega$,  let us fix a finite set of points $\{y_1, \, \ldots \, ,y_N\}$ on 
$\partial B_r(x)\cap \Omega$ such that 
$$
\partial B_r(x)\cap \Omega\subset \bigcup_{i=1}^{N} B_{\delta}(y_i),
$$
 and   let us set the function 
$w_r^x:  \Om \to \R$ defined by
\begin{equation}
w_r^x(z) := \max_i w_r^{x,y_i}(z).
\end{equation}
By construction, the following facts hold:
\vskip4pt
\begin{itemize}
\item[1)]
\vskip4pt
$\supess_{B_r(x)\cap \Omega} f(z,D w_r^x(z))  \le 1$;
\item[2)] 
\vskip4pt 
$w_r^x(z)> v(z)+\varepsilon$ for every $z\in \partial B_r(x)\cap \Omega$

(in fact,  if $z\in \partial B_r(x)\cap \Omega$ then there exists $y_i$ such that 
$z\in B_{\delta}(y_i)$. Therefore by (\ref{stima*}) 
$w_r^x(z)\geq w_r^{x,y_i}(z)> v(z)+\varepsilon);$
 \item[3)]
\vskip4pt
$w_r^x(z) < v(z)+\varepsilon$ for every $z\in  B_{\delta}(x)\cap \Om$ (it follows by  (\ref{stima**})).
\end{itemize}
\vskip4pt
Now let $Z_r$ be a finite set of points of $\Om$ such that 
$$
\Omega\subset \bigcup_{z\in Z_r} B_{\delta}(z),
$$
and consider the function $w_r:\bigcup_{z\in Z_r} B_{r}(z) \cap \Omega \to \R$ defined by
\begin{equation} \label{minpoi}
w_r(x):= \min_{z\in Z_r\cap B_r(x)} w_r^z(x).
\end{equation} 
By properties $2)$ and $3)$ above it follows that $ w_r$ is continuous on $\Omega$ (one may proceed by induction on the cardinality of $Z_r$).

Moreover, for almost every $x$ in $\Omega$,  $D w_r(x)$ coincides with $D w^z_r(x)$ for some $z\in Z_r$ and this implies that $w_r\in W^{1,\infty} (\Omega)$ and $F(w_r)\le 1$.

Now let us prove that $\|w_r -v\|_{\li} \to 0$.
To this aim, let us fix $x\in \Omega$, and let $z\in B_r(x)\cap \Omega$ be such that $w_r(x)= w_r^z(x)$.
Recalling that by construction $w_r^z(z)=v(z)$, by  using  (\ref{dFbzero}) and Lemma \ref{Lipgeod}, and by the regularity of $\partial\Omega$ , we conclude that there exists a constant $C_{\Omega}>0$  such that 
\begin{eqnarray*}
|w_r(x) - v(x)|&\le& |w_r^z(x)-w_r^z(z)|+|w_r^z(z)-v(x)|\\
&=&
|w_r^z(x)-w_r^z(z)|+|v(z)-v(x)|\\
&\le&  \max \{d_F(x,z), d_F(z,x)\}+ \max \{d_F(x,z), d_F(z,x)\}\le 2C_\Omega \beta r.
\end{eqnarray*}
Therefore, 
%for every sequence $(r_n)$ converging to $0$, 
the sequence  $v_n:= w_{1/n}$ satisfies the thesis. 
%Rimane da estendere ogni $w_r$ su tutto $\Omega$ e verificare che esiste, a meno di sottosuccessioni, una funzione $\tilde v$  limite in $\wi $ delle $w_r$ estese, che coincide con $v$ su $\Omega$.  

\qed

In the following proposition  we  provide  a characterization of functional $F$ on $\wi$ whose $d^\la_F$ satisfies 
\beq\label{alto} d^{\la}_F(x,y)\leq \beta(\lambda)|x-y|_\Omega\,\quad  \forall x,y\in\Omega.
 \eeq 
 
\begin{prop}\label{altoprop} Let $\Omega\subseteq \R^N$ be a connected  open set  and let $F:\wi\to \bar\R$ be a  functional. Let  $\lambda>\inf_{\wi}F$. Then  the following facts are equivalent:
\begin{itemize}
\item[(i)] there exists $\beta=\beta(\lambda)>0$ such that
  \eqref{alto} holds;
  \item[(ii)]   there exists $\beta=\beta(\lambda)>0$  such that \beq\label{weakcoerc}\{u\in\wi\ :\ F(u)\le \lambda\}\subseteq \{ u\in \wi \ : \ \|Du\|_{\infty}\leq \beta(\la) \}.\eeq
\end{itemize}
\end{prop}
\proof 
$(i)\Longrightarrow (ii)$ Let $u\in \wi$ be such that $F(u)\le \lambda $. Then for every $x,y\in \Omega$ 
$$|u(x)-u(y)|\leq d^{\la}_F(x,y)\vee d^{\la}_F(y,x)\leq \beta(\lambda)|x-y|_\Omega$$
and, thanks to Lemma \ref{Lipgeod}, this implies  $\|Du\|_{\infty}\leq  \beta(\lambda).$

$(ii)\Longrightarrow (i)$  Thanks to Lemma \ref{Lipgeod},  for every $u\in E_\la$, we have that  $$ u(x)-u(y)\leq ||Du||_{\infty} |x-y|_{\Omega}\leq \beta(\lambda) |x-y|_{\Omega} ;$$  
hence passing to the supremum on $u\in E_\la$ we get  
$$d^{\la}_F(x,y)\le \beta(\lambda) |x-y|_{\Omega}.
$$
\qed

\begin{rem}\label{disrelax}{\rm
Note that, under assumption  \eqref{weakcoerc} or, equivalently, under assumption  \eqref{alto},  by Remark \ref{lsc}(2)-(3) we have that \beq\label{relril}  {\Gamma_{w^*}(F)}\leq {\Gamma_{w^*_{seq}}(F)}={\Gamma_{\tau_{\infty}}(F)}\eeq  and,  taking into account \eqref{dis1},  for every $\forall \lambda >\inf_{\wi}F $ it holds 

$$ d^{\la}_{\Gamma_{w^*}(F)}(x,y)\geq d^{\la}_{\Gamma_{w^*_{seq}}(F)}(x,y)=d^{\la}_{\Gamma_{\tau_{\infty}}(F)}  (x,y)\geq d^{\la}_{F}  (x,y)\quad \forall x,y\in\Omega.
$$
We emphasize that  in general, for any $\tau \in\{ \tau_{\infty}, w^*, w^*_{seq}  \}$, the distance  $d^{\la}_{\Gamma_\tau{(F)} }$ does not coincide with $d^{\la}_{F} $ (see Example \ref{gap}).}

 \end{rem}

In the following proposition we prove that, under assumption  \eqref{weakcoerc},   the  distances $d^{\lambda}_{\Gamma_{\tau}(F)}$ can be obtained as infimum of the distances associated to $F$  for any $\tau \in \{\tau_{\infty}, w^*_{seq}\}$. In particular, when  (\ref{dFbzero}) holds, we get that  the different quotients associated to the relaxed  functional $\Gamma_{\tau}(F)$ can be obtained as supremum of the different quotients associated to $F$ (see \eqref{rapprR} for a precise statement).  

Afterwards,  in Proposition \ref{altoprop1} we will exploit \eqref{relaxdis}  in order  to  provide  a characterization of functionals $F$ on $\wi$ whose $(d^\la_F)_{\la}$ satisfies
\beq \label{basso} d^{\lambda}_{F}(x,y)\geq \alpha (\lambda) |x-y| \quad \forall x,y\in \Omega.\eeq
We underline that in all these results  a relevant hypothesis  is that  $F$ is translation invariant, that is
\beq\label{trasl}F(u+c)=F(u) \quad \forall\,u \in \wi,\, \forall\,c\in \R.\eeq
while we do not need a priori $F$ to be a supremal functional.

\begin{prop} \label{infdist} Let $\Omega\subseteq \R^N$ be a connected open set with Lipschitz continuous boundary and let $F:\wi\to \bar\R$  be a   translation invariant
 functional  satisfying \eqref{weakcoerc} for every  $\lambda >\inf_{\wi} F$. Then   
for every $\la>\inf_{\wi} F$ and for any $\tau \in \{\tau_{\infty}, w^*_{seq}\}$ we have  

\beq\label{relaxdis}\inf_{\la'>\la}d^{\lambda'}_F(x,y)=d^{\lambda}_{\Gamma_{\tau}(F)} (x,y) \qquad \forall x,y\in \Omega.
 \eeq
In addition, if  $(d^{\la}_F)_{\la}$ satisfies \eqref{dFbzero} holds, then  for every $\la>\inf_{\wi} F$ and for any $\tau \in \{\tau_{\infty}, w^*_{seq}\}$ we have  

 \beq\label{rapprR}
 R_{d^{\lambda}_{\Gamma_{\tau}(F)}}(u)=\sup_{\la'>\la} R_{d^{\lambda'}_F}(u)\qquad \forall u\in\wi.
 \eeq
\end{prop}
 \proof 
By  Remark \ref{disrelax} it holds  $\Gamma_{\tau_{\infty}}(F)=\Gamma_{w^*_{seq}}(F)$, hence   it is sufficient to show the thesis when  $\tau=\tau_{\infty}$. 
In order to prove \eqref{relaxdis}, let us  fix $\lambda >\inf_{\wi} F$. Let  $v\in \wi$  be such that $\Gamma_{\tau _{\infty}}(F)(v)\leq \la$. 
 Then for $\la'>\la$ there exists $(v_n)_n\subseteq \wi$ uniformly converging to $v$ such that $F(v_n)\leq \la'$ for every $n\in \N$. Thus, for all $x,y\in \Omega$, passing to the limit in the inequality 
$$d^{\lambda'}_F(x,y)\geq v_n(x)-v_n(y)\  \ \forall n\in \N, \forall x,y\in \Omega$$ 
we have 
$$d^{\lambda'}_F(x,y)\geq v(x)-v(y) \quad \forall x,y\in \Omega.$$
Hence $$d^{\lambda'}_F(x,y)\geq \sup \{ v(x)-v(y) :\ \Gamma_{\tau_{\infty}}(F)(v)\leq \lambda\}= d^{\lambda}_{\Gamma_{\tau_{\infty}}(F)}(x,y) $$
and passing to the infimum, it follows
$$\inf_{\la'>\la}d^{\lambda'}_F(x,y)\geq d^{\lambda}_{\Gamma_{\tau_{\infty}}(F)} (x,y) \quad \forall x,y\in \Omega.$$
Vice versa, let us fix $x,y\in \Omega$ and let $(\la_n)$ be a decreasing sequence converging to $\la$ such that $$\lim_{n\to \infty} d^{\lambda_n}_F(x,y)=\inf_{\la'>\la}d^{\lambda'}_F(x,y).$$ 
Let us choose $u_n\in \wi$  such that $F(u_n)\leq \la_n$ with $d^{\lambda_n}_F(x,y)\leq u_n(x)-u_n(y)+\frac 1 n$. 
Thanks to \eqref{trasl} we can assume  that $\int_{\Om} u_n dx=0$ for every $n\in \N$.
Since the sequence $(u_n)\subseteq \wi$ has bounded gradients, then, up to subsequence,  $(u_n)$ uniformly converges to a function $u_0\in \wi$. By lower semicontinuity we have  $\Gamma_{\tau_{\infty}}(F)(u_0)\leq \la$ and   $$ \inf_{\la'>\la}d^{\lambda'}_F(x,y) =\lim_{n\to +\infty} u_n(x)-u_n(y)= u_0(x)-u_0(y)\leq  d^{\lambda}_{\Gamma_{\tau_{\infty}}(F)}(x,y).$$
In order  to show \eqref{rapprR} we note that,   if  $\la>\inf_{\wi}F,$ thanks to (\ref{dis1}) and (\ref{dFbzero}), we have that 
$$ \alpha (\lambda) |x-y| \leq d^{\lambda}_{F}(x,y)\ \leq d^\lambda_{\Gamma_{\tau}(F)}(x,y).$$
Therefore $$d^\lambda_{\Gamma_{\tau}(F)}(x,y) =0  \Longleftrightarrow x=y$$ and thus, for every $u \in \wi$, we have that  

\begin{eqnarray*}
 R_{d^{\lambda}_{\Gamma_{\tau}(F)}}(u)&=&\sup_{x,y\in \Omega,\,x\ne y}\frac{u(x)-u(y)}{d^\lambda_{\Gamma_{\tau}(F)}(x,y)}\\&=& \sup_{x,y\in \Omega,\, x\ne y}\frac{u(x)-u(y)}{ \inf\limits_{\la'>\la}d^{\lambda'}_F(x,y)}\\
 &=&  \sup_{x,y\in \Omega, \, x\ne y}\sup\limits_{\la'>\la} \frac{u(x)-u(y)}{ d^{\lambda'}_F(x,y)}\\
 &=& \sup\limits_{\la'>\la} \sup_{x,y\in \Omega, \, x\ne y}\frac{u(x)-u(y)}{ d^{\lambda'}_F(x,y)}
 = \sup_{\la'>\la} R_{d^{\lambda'}_F}(u). 
\end{eqnarray*}

\qed

%Proposition \ref{altoprop} and \eqref{relaxdis} in Proposition \ref{infdist},} we are in position to characterize supremal  functionals $F$ on $\wi$ whose associated $d^\la_F$ satisfy (\ref{dFbzero}). 
%In the following proposition  we  provide  a characterization of functional $F$ on $\wi$ whose $d^\la_F$ satisfies 

\begin{prop}\label{altoprop1} Let  $\Omega$ be a bounded connected open set with Lipschitz continuous boundary  and let $F:\wi\to \bar\R$ be a   translation invariant
 functional  satisfying \eqref{weakcoerc} for every $\la>\inf_{\wi} F$.
If 
$$ \lim_{\xi\to 0} {\Gamma_{\tau}(F)} (u_{\xi})= \inf_{\wi} F $$
 where  $u_{\xi}(x):=\xi\cdot x$ and  $\tau \in \{\tau_{\infty}, w^*_{seq}\}$,  then 
for every $\la>\inf_{\wi} F$ there exists $\alpha (\lambda)>0$ such that  \eqref{basso} holds.
\end{prop}
\proof 
%$(i)\Longrightarrow (ii)$ Let $u\in \wi$ be such that $F(u)\le \lambda $. Then for every $x,y\in \Omega$ 
%$$|u(x)-u(y)|\leq d^{\la}_F(x,y)\vee d^{\la}_F(y,x)\leq \beta(\lambda)|x-y|_\Omega$$
%and, thanks to Lemma \ref{Lipgeod}, this implies  $\|Du\|_{\infty}\leq  \beta(\lambda).$
%Therefore it is sufficient to choose $r= \beta(\lambda)$.
%
%
%$(ii)\Longrightarrow (i)$ Without loss of generality, we can suppose that  there exists a continuous,  increasing and one-to-one function $\phi:[0,+\infty)\to [0,+\infty)$ such that $$\lim_{t\to +\infty} \phi(t)=+\infty$$ and 
%$$F(u)\geq \|\phi(Du)\|_{\infty}.$$ 
% Note that by hypothesis we deduce  $d^{\la}_F(x,y)<+\infty$. 
% Let $\epsilon>0$ and let  $u_0\in \wi$ be such that  $F(u_0)\le \lambda $ and 
%$ d^{\la}_F(x,y) \leq u_0(x)-u_0(y)+\epsilon.$
%   Thanks to the coercivity assumption we have that $$\phi(\|Du_0\|_{\infty})= \|\phi(Du_0)\|_{\infty} \leq F(u_0)\leq \lambda$$ that implies 
%$$\|Du_0\|_{\infty}\leq \phi^{-1}(\lambda).$$
%So we obtain that 
%$$d^{\la}_F(x,y)\leq u_0(x)-u_0(y)+\epsilon\leq \|Du_0\|_{\infty} |x-y|_{\Omega}+\epsilon\leq   \phi^{-1}(\lambda) |x-y|_{\Omega}+\epsilon.$$\
%When $\epsilon \to 0$ we obtain (\ref{alto}) with $\beta=\phi^{-1}(\lambda)$. 
%
%

Let  $\la>\inf_{\wi} F$ be fixed and choose any value $\tilde \lambda\in \R$ such that $\inf_{\wi} F<\tilde \lambda<\la $. By the definition of limit there exists  $\delta>0$, depending only on $\tilde \lambda$,  such that for every $\xi\in \R^N$ such that $|\xi|\leq \delta$ it holds $\Gamma_{\tau_{\infty}}(F)(u_{\xi})\leq \tilde \lambda $.  Set $\alpha=\alpha(\tilde \lambda)=\delta/2$. We claim that for every $x,y\in \Omega, x\ne y$
$$d^{\lambda}_F(x,y)\geq\alpha  |x-y|.$$
Indeed, let $x,y\in \Omega, x\ne y$ and define $\eta:=\alpha \frac{(x-y)}{|x-y|}$. Since $|\eta|= \alpha<\delta$ we get 
$$\Gamma_{\tau_{\infty}}(F)(u_{\eta})\leq  \tilde\lambda.$$
By definition of $d^\la_{\Gamma_{\tau_{\infty}}(F)}$ and by (\ref{relaxdis}) this implies that 
  $$d^{\lambda}_F(x,y)\ge\inf_{\gamma>\tilde\la}d^{\gamma}_F(x,y)=d^{\tilde\la}_{\Gamma_{\tau_{\infty}}(F)}(x,y)\geq u_{\eta}(x)-u_{\eta}(y)=\alpha  |x-y|. $$\qed

Finally we can give a characterization of supremal functionals $F$ whose  distances $d^\la_F$ satisfy \eqref{dFbzero}. 

\begin{thm} \label{characterization} Let  $\Omega$ be a bounded connected open set with Lipschitz continuous boundary and let  $F$ be a supremal functional of the form \eqref{sfzero}. Then the following facts are equivalent

\begin{itemize}
\item[(i)]$(d^{\la}_F)_{\la}$ satisfies  \eqref{dFbzero};
\item[(ii)] $F$ satisfies  \eqref{weakcoerc}  $\forall \la>\inf_{\wi} F$ and  
\beq\label{lim0} \lim_{\xi\to 0} {\Gamma_{\tau}(F)} (u_{\xi})= \inf_{\wi} F= {\Gamma_{\tau}(F)}(0) \quad \forall  \tau \in \{  \tau_{\infty}, w^*_{seq} \} .\eeq 
\end{itemize}
\end{thm}

\proof 

Assume that  $(d^{\la}_F)_{\la}$ satisfies  (\ref{dFbzero}).  By Proposition \ref{altoprop}  we get that $F$ satisfies \eqref{weakcoerc}  $\forall \la>\inf_{\wi} F$. Moreover by Remark \ref{disrelax}  it is sufficient to show \eqref{lim0} with $\tau=\tau_{\infty}$.
Fix $\la>\inf_{\wi} F$.
Since 
  $$R_{d^{\lambda}_F }(u_{\xi})=\sup_{x,y\in \Omega,\,x\ne y}\frac{u_{\xi}(x)-u_{\xi}(y)}{d^\lambda_{F}(x,y)}\leq |\xi| \sup_{x,y\in \Omega,\,x\ne y}\frac{|x-y|}{d^\lambda_{F}(x,y)}\leq  |\xi| \sup_{x,y\in \Omega,\,x\ne y}\frac{|x-y|}{\alpha (\lambda) |x-y|}=\frac{|\xi|}{\alpha (\lambda) }, $$
if $|\xi|<\alpha (\lambda)$ then, by Proposition  \ref{keylemma},  there exists $(u_n^\xi)_n\subseteq \wi$ such that $u_n^\xi\to u_{\xi}$ uniformly and $F(u_n^\xi)\leq \la $ for any $n\in \N$. 
 This implies     $$   {\Gamma_{\tau_{\infty}}(F)} (u_{\xi})\leq \liminf_{n\to \infty }F(u_n^\xi)\leq \la$$
for every  $|\xi|<\alpha (\lambda)$.
Up to passing to the limsup  as $\xi \to 0$  we deduce that  
$$
\limsup_{\xi\to 0} {\Gamma_{\tau_{\infty}}(F)} (u_{\xi})\leq\la. 
$$
Eventually,  we let $\la\to \inf_{\wi} F$  and we get that
  $$ \inf_{\wi} F\leq {\Gamma_{\tau_{\infty}}(F)} (0) \leq \liminf_{\xi\to 0} {\Gamma_{\tau_{\infty}}(F)} (u_{\xi})
  \leq \limsup_{\xi\to 0} {\Gamma_{\tau_{\infty}}(F)} (u_{\xi})\leq \inf_{\wi} F.$$

The converse implication follows by Propositions \ref{altoprop} and \ref{altoprop1}.
   \qed
	
\begin{rem}\label{alto2} 
{\rm We remark that hypotheses \eqref{alto}, \eqref{basso} can be easily infered for functionals $F$ with linear growth from above and below respectively.  Indeed,  
 if there exists $\beta>0$  such that $$F(u)\geq \beta \|Du\|_{\li(\Omega)}\quad \forall\, u\in \wi$$ then 
for every $\lambda>0$ and for every $x,y\in \Omega$ it holds
$$ d^{\lambda}_{F}(x,y)\leq \frac \lambda {\beta}  |x-y|_\Om. $$ 
On the other hand if there exists $\alpha>0$ such that  
$$
F(u)\leq \alpha\|Du\|_{\li(\Omega)}\quad \forall\, u\in \wi
$$
 then, for every $\lambda>0$,  $\alpha(\la)$ in \eqref{basso} can be chosen equal to $ \displaystyle \frac \lambda {\alpha}$.
% Indeed,   it is sufficient to choose $\eta:=\frac \lambda C \frac{(x-y)}{|x-y|}$ above. Since  $F(u_{\eta})\leq \lambda$ we obtain  that for every $x,y\in \Omega$ 
%  $$d^\la_{F}(x,y)\geq u_{\eta}(x)-u_{\eta}(y)=\frac \lambda C  |x-y|. $$
}
\end{rem}
\bigskip

We conclude this section by giving  an alternative representation  for $d^{\la}_F$. To this aim we introduce for any subset $A$ of a topological vector space $(X, \tau)$ 
\beq\label{tauco}
\tau\hbox{-}\conv(A):=\bigcap \{B \hbox{ $\tau$-closed and convex subset of }X\, :\, B\supseteq A\}. 
\eeq

\begin{prop} \label{chiusura} Let $\Omega\subseteq \R^N$ be a connected open set with Lipschitz continuous boundary. Let $F:\wi\to \bar\R$ be a  functional and let $\tau$ be one of  the topologies $\tau_{\infty}, w^*, w^*_{seq}$. Then $\forall \lambda >\inf_{\wi} F$, $\forall x,y\in\Omega$  it holds 
\begin{equation}\label{distconv}
d^{\la}_F(x,y)=\sup\{u(x)-u(y):\, u\in \tau\hbox{-}\conv(\{v\in \wi: F(v)\le \lambda\})\} , 
\end{equation}
Moreover, if  $F$ is a supremal functional and  $ (d^\la_F)_{\la}$ satisfies (\ref{dFbzero}) then 
\beq\label{conlevelset} \tau\hbox{-}\conv\Big(\{ u\in\wi : \  F(u)\leq\lambda  \}\Big)
=\{u\in\wi : \  R_{d^{\lambda}_F}(u)\leq 1  \}. \eeq
\end{prop} 
\proof
Note that, thanks to the Banach-Dieudonne-Krein-Smulian Theorem (see \cite[Theorem 3.33]{Brezis}) and the definition of $\tau\hbox{-}\conv (A)$,    $w^*\hbox{-}\conv(A)=w^*_{seq}\hbox{-}\conv(A)$ for any given set $A\subseteq \wi$. 
Hence it suffices to show \eqref{distconv} and \eqref{conlevelset} for $\tau \in \{\tau_{\infty}, w^*_{seq}\}$. 

Moreover, \eqref{distconv} follows straigthforward once we prove that  
\begin{equation}\label{coA}\sup\{u(x)-u(y):\, u\in A\} =\sup\{u(x)-u(y):\, u\in \tau\hbox{-}\conv(A)\} \end{equation} 
for any given set $A\subseteq \wi$.  
To this aim let us define 
$$\conv (A):=\left\{\sum_{i=1}^m\lambda_i u_i\, : \lambda_i\ge 0, u_i\in A, \sum_{i=1}^m \lambda_i=1, m\in \N\right\}.
$$

It can be easily verified that $\conv(A)$ is a convex set containing $A$ and the $\tau$-closure of $\conv (A)$  coincides with $\tau\hbox{-}\conv(A)$ for $\tau \in \{\tau_{\infty}, w^*_{seq}\}$.  For such topologies one may argue by sequences  and prove easily that  
$$\sup\{u(x)-u(y):\, u\in \conv (A)\} =\sup\{u(x)-u(y):\, u\in \tau\hbox{-}\conv(A)\}. $$
Thus it suffices to show that 
$$\sup\{u(x)-u(y):\, u\in A\}=\sup\{u(x)-u(y):\, u\in \conv (A)\}.$$
One inequality follows by comparison. To get the converse inequality, it is enough to note that for any convex combination  $\sum_{i=1}^m\lambda_i u_i\in \conv (A)$ and for any $ x,y\in\Omega$ we have  
\begin{eqnarray*}
\sum_{i=1}^m\lambda_i u_i(x)-\sum_{i=1}^m\lambda_i u_i(y)&=& \sum_{i=1}^m\lambda_i (u_i(x)-u_i(y))\\
&\le& 
\sum_{i=1}^m\lambda_i \sup\{u(x)-u(y):\, u\in A\}=\sup\{u(x)-u(y):\, u\in A\}. 
\end{eqnarray*}

It remains to show \eqref{conlevelset}. To this aim, let  $\tau \in \{\tau_{\infty}, w^*_{seq}\}$ and let $u\in\wi$ be such that  $F(u)\leq\lambda $. Then, by definition of $d^{\lambda}_F$, it holds $R_{d^{\lambda}_F}(u)\leq 1 $.
This in turn implies that $\{ u\in\wi : \  F(u)\leq\lambda  \}\subseteq \{u\in\wi : \  R_{d^{\lambda}_F}(u)\leq 1  \}$ and, since $R_{d^{\lambda}_F}$ is a convex functional, we get that $$\tau\hbox{-}\conv\Big(\{ u\in\wi : \  F(u)\leq\lambda  \}\Big)
\subseteq \{u\in\wi : \  R_{d^{\lambda}_F}(u)\leq 1  \}. 
$$
\noindent On the other hand, if  $R_{d^{\lambda}_F}(v)\leq 1$ then, for $0<\epsilon<1$,  $R_{d^{\lambda}_F}(\epsilon v)=\epsilon R_{d^{\lambda}_F}(v)< 1$ and, by Proposition \ref{keylemma}, we have that
there exists a sequence $\{v_n\}\subset\wi$ converging to $\epsilon v$ in $\li$ with $F(v_n) \le \la$.  Then $\epsilon v\in {\tau} \hbox{-}\conv\Big(\{ u\in\wi : \  F(u)\leq\lambda  \}\Big)  $ that implies, when $\epsilon\to 1$,  that $v\in {\tau}\hbox{-}\conv\Big(\{ u\in\wi : \  F(u)\leq\lambda  \}\Big)  $.

\qed

\section{The proofs}\label{proofs}
\noindent {\sl Proof of Theorem \ref{relax-1}.}  We first prove \eqref {rapplevels} for  $\tau=\tau_{\infty}$.
Let $\lambda >\inf_{\wi} F=\inf_{\wi } \Gamma_{\tau}(F)$ be fixed.  

We note  that, by definition, 
\beq \label{incl1}
 \{ u\in\wi : \  \Gamma_{\tau} (F)(u)\leq\lambda  \}\subseteq \{ u\in W^{1,\infty}(\Omega)\ : \ R_{d^{\lambda}_{\Gamma_{\tau}(F)}}(u)\leq 1   \}.
\eeq
Moreover, by \eqref{rapprR}  in Proposition \ref{infdist}  it easily follows that   \beq \label{incl2}\{ u\in W^{1,\infty}(\Omega)\ : \ R_{d^{\lambda}_{\Gamma_{\tau}(F)}}(u)\leq 1  \}= \{ u\in W^{1,\infty}(\Omega)\ : \ \sup_{\la'>\la} R_{d^{\lambda'}_F}(u)\leq 1   \}.\eeq

If we prove that  
\beq\label{incl3}
 \{ u\in W^{1,\infty}(\Omega)\ :\ \sup_{\la'>\la} R_{d^{\lambda'}_F}(u) \leq 1 \}\subseteq \{ u\in\wi : \  \Gamma_{\tau} (F)(u)\leq\lambda  \}
\eeq
 then, by \eqref{incl1}, (\ref{incl2}) and (\ref{incl3})  we get  (\ref{rapplevels}).

Let $u\in \wi$ be such that $R_{d_F^{\lambda'}}(u)\leq 1$ for every $\lambda'>\lambda$.  Then,
 for every fixed $0<\theta<1$ we have that  $$R_{d_F^{\lambda'}}(\theta u)=\theta R_{d_F^{\lambda'}}(u) \leq \theta <1.$$	 With fixed $\lambda '>\lambda$ we  apply Proposition \ref{keylemma}  and get that   there exists a sequence $(u^{\lambda',\theta}_n)_n\subset W^{1,\infty}(\Omega)$ converging to $u$ in $L^{\infty}(\Omega)$ such that  $F(\theta u^{\lambda',\theta}_n) \le \lambda' $. Then $\Gamma_{\tau_{\infty}}(F)( \theta u) \le \lambda'$ for every $\lambda' >\lambda$ and for every $0<\theta<1$ which implies 
$ \Gamma_\tau (F)(\theta  u) \le  \lambda.$ Letting $\theta\to 1^-$, by lower semicontinuity, we get $ \Gamma_\tau (F)( u) \le  \lambda.$

Now we claim that the functional $\Gamma_{\tau_{\infty}} (F)$ is level convex. 

Indeed, since  $R_{d_{\Gamma_{\tau} (F)}^{\lambda}}$ is convex, we conclude that the sublevel sets of   $\Gamma_{\tau}(F)$ are convex for $\lambda >\inf_{\wi} F$. For $\lambda =\inf_{\wi} F$ it suffices to note that  
$$\left\{u\in\wi : \Gamma_{\tau}(F)(u)= \inf_{\wi} F\right\}=\bigcap_{\lambda >\inf\limits_{\wi} F}\left\{u\in \wi: \Gamma_\tau (F)(u)\le \lambda\right\}.$$ 
Finally,  thanks to Proposition \ref{altoprop}, the functional $F$ satisfies (\ref{weakcoerc}) and, by applying  Remark \ref{lsc}(2),  we have that $\Gamma_{w_{seq}^*}(F) \equiv \Gamma_{\tau_{\infty}} (F)$. Therefore $ \Gamma_{w_{seq}^*}(F)$  is level convex and  by Proposition \ref{weakseq} we can conclude that $$\Gamma_{\tau_{\infty}} (F)\equiv \Gamma_{w^*}(F)\equiv \Gamma_{w_{seq}^*}(F).$$  
Hence \eqref {rapplevels} holds for any  $\tau\in \{\tau_{\infty}, w^*, w^*_{seq}\}$.
  
%Finally,  \eqref{rapplevels2} easily follows from \eqref{rapplevels} when  $F$ is  $\tau$-lower semicontinuous.
\qed

\noindent{\sl Proof of Theorem \ref{generGPP}.}  First of all we note that given a functional $G$  positively $1$-homogeneous  for any strictly positive  value $\lambda >\inf_{\wi} G$ and $\forall x,y\in \Omega$ it holds 
\beq\label{homogen} d^\lambda_G(x,y)=\sup\{u(x)-u(y) : G(u)\le \lambda\}=\lambda \sup\left\{\frac{u}{\lambda}(x)-\frac{u}{\lambda}(y) : G(\frac{u}{\lambda})\le 1\right\}=\lambda d^1_G(x,y) .
\eeq
In particular the scaling  in \eqref{homogen} can be applied for $\lambda >0$ both to $F$ and $\Gamma_\tau(F)$, as it can be easily verified that $\Gamma_\tau(F)$ is also positively $1$-homogeneous. 

\noindent We now claim that a supremal functional $F$ satisfying hypotheses $(1), (2)$ is non negative. In order to prove this claim  we assume in the sequel $\tau=\tau_\infty$. 
%In order to prove this claim  we first observe that by Theorem \ref{relax-1}  the relaxed functional $\Gamma_\tau(F)$  does not depend on the choice of the topology $\tau$ so we may choose $\tau=\tau_\infty$.  

Hence we may apply Proposition \ref{infdist} and obtain by \eqref{relaxdis} the equality 
\beq\label{uguadist}d^1_{\Gamma_{\tau}(F)} (x,y)=\inf_{\la>1}d^{\lambda}_F(x,y)= \inf_{\la>1}\lambda \, d^1_F(x,y)= d^1_F(x,y)\quad \forall x,y\in \Omega.
 \eeq

\noindent By using hypothesis $(2)$ on $d^1_F=d^1_{\Gamma_{\tau}(F)} $,  by Proposition \ref{altoprop},  we infer that the sublevel set $\{u\in \wi: \Gamma_\tau(F)(u)\le 1\}$ consists of functions with equibounded gradients.  Moreover this property is inherited by inclusion by any sublevel set $\{u\in \wi: \Gamma_\tau(F)(u)\le \lambda\}$ for $\lambda <1$. 

Let us assume by contradiction that $\inf_{\wi}F<0$. By \eqref{lim0} in Proposition \ref{characterization} we get that 
$$\lim_{\xi\to 0}\Gamma_\tau(F)(u_{\xi})=\Gamma_\tau(F)(0)=\inf_{\wi}F<0 .$$ Hence, there exists $\bar\xi\ne 0$ such that  $\Gamma_\tau(F)(u_{\bar\xi})=\eta <0$. By the positively $1$-homogeneity we also infer that $\forall \mu >1$ 
$\Gamma_\tau(F)(u_{\mu\bar\xi})=\Gamma_\tau(F)(\mu u_{\bar\xi})=\mu\eta <\eta$ and this leads to  a contradiction since the set 
$\{u_{\mu\bar\xi},\, \mu>1\}$ has non uniformly bounded gradients and thus it cannot be contained in the  sublevel set $\{u\in\wi : \Gamma_\tau(F)(u)\le \eta\}$. 

Hence the claim is proved and by \eqref{homogen} we have $d^{\lambda}_F(x,y)=\lambda d^1_F(x,y)$ $\forall x,y\in \Omega$ and $\lambda >\inf_{\wi}F$ as  $\inf_{\wi}F\ge 0$. Therefore $F$ satisfies \eqref{dFbzero} and  by the first equality of \eqref{rapplevels} in Theorem \ref{relax-1}   we get 
\beq\label{livelli}
\{ u\in\wi : \  \Gamma_{\tau} (F)(u)\leq\lambda  \} =\{u\in\wi : R_{d^1_F}(u)\leq \lambda\} 
\eeq 
for $\lambda >\inf_{\wi}F\ge 0$. 

Thus the statement $\Gamma_\tau (F) =R_{d^1_{\Gamma_{\tau} (F)}}$ is proved once  we show that 
$$\inf_{\wi} F= \inf_{\wi}\Gamma_\tau (F) =\inf_{\wi}R_{d^1_{\Gamma_{\tau} (F)}}.
$$ 
To this aim we actually prove that the infima above are minima and they are achieved on constant functions, so that  $\min_{\wi}\Gamma_\tau (F) =0=\min_{\wi}R_{d^1_{\Gamma_{\tau} (F)}}$.   

In order to get  $\inf _{\wi}R_{d^1_{\Gamma_{\tau} (F)}}=0$ it suffices to notice that, thanks to hypothesis $(2)$, the values  $d^1_F(x,y)$ are finite $\forall x,y\in \Omega$, thus, up to exchanging the role of $x$ and $y$,  the functional $R_{d^1_{\Gamma_{\tau} (F)}}$ is non negative; moreover $R_{d^1_{\Gamma_{\tau} (F)}}(u)=0$  on any constant function $u$. It remains to prove that $\inf _{\wi} \Gamma_\tau (F)= 0$.  This can be obtained by using \eqref{livelli} with $\lambda =1$  and noticing that the constant function $\bar u=0$ belongs to the sublevel set $ \{R_{d^1_F}\leq 1\}$. Hence  $\Gamma_\tau (F)(\bar u)\le 1$ and, by the positively $1$-homogeneity,  it also holds 
$$\Gamma_\tau (F)(0)=\Gamma_\tau (F)(2 \bar u)= 2\Gamma_\tau (F)(\bar u)=2 \Gamma_\tau (F)(0)$$
and this implies $\Gamma_\tau (F)(0)=0$.  
\qed

\bigskip

In order to give the proof of Theorem \ref{relax0} we need to show some preliminary results. The first one concerns an approximation  via $\Gamma$-convergence of  the $\tau$-lower semicontinuous envelope of a   non negative  functional $F$ on $\wi$ through a sequence of coercive  functionals. 

\begin{prop}\label{approxrelax}  Let $\om$ be   an open subset of $\R^N$.  Let  $F:\wi \times \R^N\to [0,+\infty]$ be a functional 
 and for every  $n\in \N$ let $F_n:\wi\to [0,+\infty]$ be the functional   defined  by 
\beq\label{Gn}
F_{n}(u):=F(u)\vee \frac{1}{n} ||Du||_{L^{\infty}(\om)}.
\eeq
 Then
\begin{itemize} 

\item [\rm (1)] the sequence   $ (F_n)_{n\in \N}$ is decreasing and pointwise converges to $F$;
\item [\rm (2)] $F$ is a level convex functional if and only if  $F_n$ is  level convex  for every $n\in\N$;
\item  [\rm (3)] for any topology $\tau\in\{\tau_\infty, w^*, w^*_{\rm seq}\}$ the sequences  $( F_n)_{n\in\N}$ and $(\Gamma_{\tau}( F_n))_{n\in\N}$ $\Gamma(\tau)$-converge to $\Gamma_{\tau}( F)$ in $\wi$.
\end{itemize}

\end{prop} 

\proof {\rm (1)} derives from the fact that $F\geq 0$. In order to prove {\rm (2)} we first note that if  $F$ is level convex then, for every $n\in \N$,  then the sublevel set $$\{ u\in \wi: \ F_n(u)\leq \la\}= \{ u\in \wi: \ F(u)\leq \la\}\cap \{ u\in \wi: \ ||Du||_{\infty}\leq \la\}$$ is convex. Proposition \ref{lcGammalimite} yields the converse implication. Finally, {\rm (3)} follows by statements (2) and (3) in Proposition \ref{gammaprop}. 

\qed\\

 Thanks to the next proposition, if $F$ is a supremal functional then  $F_n$ defined by (\ref{Gn}) is itself a supremal functional. Note that in general the functional $F\wedge G$ is not a supremal functional, as shown in Section \ref{examples}.
 
\begin{prop}\label{specialcase}  Let $F,G:\wi\to \bar\R$ be  supremal functionals represented by the supremands $f,g:\Omega \times \R^N\to \bar\R$, respectively.  Then the functional  $F\vee G$  defined by $F\vee G(u):=F(u)\vee G(u)$ is still a supremal functional, represented by the supremand $f\vee g$.
\end{prop}
\proof Set  $H(u):=\supess_{\Omega}(f\vee g)(x,Du(x)) $ for every $u\in \wi$, we show  that  $H= F\vee G$.   The inequality  $H\geq F\vee G$ follows by definition.  In order to prove the converse inequality, let  $u\in \wi$ and $\delta>0$ be fixed and select $B_{\delta}\subset \Om$ with   $|B_{\delta}|>0$  such that  $(f\vee g)(x,Du(x))\geq H(u)-\delta$  for every $x\in B_{\delta}.$ 
Set $B^f_{\delta}:=\{x\in B_{\delta} :\,(f\vee g)(x, Du(x))= f(x, Du(x)) \}$ and  $B^g_{\delta}:=B_{\delta} \setminus B^f_{\delta}.$
If $|B^f_{\delta}|>0$, then  $F(u)\geq H(u)-\delta$ while if $|B^-_{\delta}|>0$ then  $G(u)\geq H(u)-\delta.$
In both cases $(F\vee G)(u) \geq H(u)-\delta$ for every $\delta >0$ and this entails $H\leq (F\vee G)$. 
In particular,  it follows that if $F$ is a supremal functional,  the functional $F_n$ given by (\ref{Gn}) is still a supremal functional.   Moreover if  $F(u)\geq 0$ for every $u\in \wi$ then   $F(u)=F(u)\vee 0=\supess_{\Omega} f^+(x,Du(x))$ for every $u\in\wi$ where $f^+=f\vee 0$. Therefore if $F$ is non negative we can suppose, without loss of generality,   that its supremand $f$ is non negative. \qed
\bigskip

We are in position to show Theorem \ref{relax0}. \\

\noindent {\sl Proof of Theorem \ref{relax0}.} We divide the proof into four steps. 

\noindent {\bf Step 1.}  We  assume that  $u_n\equiv 0$ for any $n\in\N$, that is $0$ is a minimum point for $F$,  and that  $F$ is coercive, i.e. $\exists \beta>0$ such that  
\beq\label{coerc}
F(u)\geq \beta ||Du||_{L^{\infty}(\Om)}\hbox { for every }u\in \wi.
\eeq
Note that, in this case, $$\Gamma_{\tau}(F)(u)\geq\min_{u\in \wi} F(u)= F(0)$$ for every $u\in \wi$ and for any $\tau$-topology.  Thus  \eqref{accaomega} implies that for any $\tau$-topology \beq\label{uxi}\min_{\wi} F  =  \Gamma_{\tau} (F)(0)\leq \liminf _{\xi\to 0} \Gamma_{\tau}(F)(u_{\xi})\leq \limsup _{\xi\to 0} \Gamma_{\tau}(F)(u_{\xi})\leq\limsup_{\xi\to 0} F (u_{\xi})=\min_{\wi} F .\eeq
\noindent By \eqref{coerc} and  \eqref{uxi}, taking into account Theorem \ref{characterization},  we deduce that   for every $\lambda >\inf_{\wi}F$  there exist $\beta (\lambda), \alpha (\lambda)>0  $ such that 
$$
 \alpha (\lambda) |x-y|  \le d^{\lambda}_{F}(x,y)\leq \beta (\lambda) |x-y|_\Omega \quad \forall x,y\in \Omega. 
$$
By applying Theorem \ref{relax-1}  we may conclude that  $\Gamma_{\tau_{\infty}}(F)$ is a level convex functional and  $$\Gamma_{\tau_{\infty}} (F)\equiv \Gamma_{w^*}(F)\equiv \Gamma_{w_{seq}^*}(F).$$   %By  Remark \ref{lsc}  it follows that $\Gamma_{\tau_{\infty}} (F)\equiv \Gamma_{w^*}(F)\equiv \Gamma_{w_{seq}^*}(F)$. Therefore it is sufficient to show that $\Gamma_{\tau_{\infty}}(F)$ is level convex. 

\noindent {\bf  Step 2.}  Here we assume the same hypotheses of Step 1 dropping  the coercivity assumption (\ref{coerc}), that is, we assume $u_n\equiv 0$ for any $n\in\N$ and $F(0)= \min_{\wi} F$. 
Without loss of generality, it is not restrictive to  suppose $F(0)=0$.  
Then we may approximate $F$ through   the sequence $(F_n)_n$ given by  $$F_n(u)= F(u)\vee  \frac 1 n  ||Du||_{\infty}.$$
 By Proposition  \ref{specialcase} $F_n$ is still  a supremal functional (represented by $f_n(x,\xi):=f(x,\xi)\vee \frac 1 n |\xi|$) and satisfies \eqref{accaomega}. By  Step 1 we get that $\Gamma_{\tau}(F_n)$ is  level convex for every $n\in \N$ when  $\tau$ is one of  the topologies $\tau_{\infty}, w^*, w^*_{seq}$.   By taking into account Proposition \ref{approxrelax}(3),  the sequence $(\Gamma_{\tau}(F_n))_n$    $\Gamma(\tau)$-converges to $\Gamma_\tau(F)$ with respect to any  topology $\tau$. Therefore, by Proposition \ref{lcGammalimite}(2),  it follows that  $\Gamma_{\tau}(F) $ is a  level convex functional for every $\tau\in \{ \tau_{\infty}, w^*, w^*_{seq}\}$.

\noindent{\bf Step 3.} We consider the general case and we  assume only the additional hypothesis   $$\inf_{u\in \wi}F(u) \in \R.$$ By taking  \eqref{accaomega} into account this implies that $F(u_n), \limsup_{\xi\to 0}  F(u_{\xi}+u_n)$ are finite for $n$ large enough. 
For every $n\in \N$ 
set 
$$c_n:= F(u_n)\vee \limsup_{\xi\to 0} F(u_{\xi}+u_n) $$ 
%where $u_n$ satisfies \eqref{accaomega} 
and define $G_n:\wi\to [0,+\infty]$ as $$G_n(u):=F(u+u_n)\vee c_n-c_n.$$
By Proposition \ref{relaxsomma} (3),(4),(5) we have that $$\Gamma_{\tau}(G_n)(\cdot)=\Gamma_{\tau}(F)(\cdot+u_n)\vee c_n-c_n.$$
Note that for every $n\in \N$  $G_n$ is a supremal functional satisfying  the assumptions of Step 2 since  $G_n(0)=0=\min\limits_{\wi} G(u)$ and 
$$\limsup_{\xi\to 0}  G_n(u_{\xi})=\limsup_{\xi\to 0}  F(u_{\xi}+u_n)\vee c_n-c_n=0 .$$
Hence $\Gamma_{\tau}(G_n)$ is level convex for any $\tau= \tau_{\infty}, w^*, w^*_{seq}$.
It easily follows that also the functional $\Gamma_{\tau}(F)\vee c_n$ is level convex for every $n\in \N$.
By passing to the pointwise limit when $n\to \infty$, by Proposition \ref{lcGammalimite} (1) we  get that $\Gamma_{\tau}(F) $ is level convex. 

\noindent{\bf Step 4.} Finally we consider the case when $\inf\limits_{u\in \wi}F(u) =-\infty$. For every $m\in \N$ let  $G_m:\wi\to \R\cup\{+\infty\}$ be given by 
$$G_m(u):=F(u)\vee (-m).$$ Note that $G_m\not\equiv +\infty$ and $\inf\limits_{\wi} G_m\in \R$. 
It is easy to show that $(u_n)_n$ is a minimizing sequence for $G_m$  and $\lim_{\xi\to 0} G_m(u_n+u_{\xi})=G_m(u_n).$ 
Since  $\inf\limits_{\wi} G_m \in \R$, by Step 3 it follows that $\Gamma_{\tau}(G_m)$ is level convex when  $\tau$ is one of  the topologies $\tau_{\infty}, w^*, w^*_{seq}$. By Proposition \ref{relaxsomma}(4) we have that $$\Gamma_{\tau}(G_m)=\Gamma_{\tau}(F\vee (-m))=\Gamma_{\tau}(F)\vee (-m).$$
By passing to the limit when $m\to \infty$ and by applying  Proposition \ref{lcGammalimite}(2) it follows that  $\Gamma_{\tau}(F) $ is a  level convex functional when $\tau =\tau_{\infty}, w^*, w^*_{seq}$. 

\qed

\noindent {\sl Proof of Corollary \ref{stability}.} Thanks to the hypotheses, by applying Theorem \ref{relax-1} or \ref{relax0}, we get  that, for any $n\in \N$ the functional $\Gamma (\tau)(F_n)$ is level convex. The thesis follows taking into account  Propositions \ref{gammaprop} (2) and \ref{lcGammalimite} (2). 
\qed

\section{Further results.}\label{repreresult}
In this section we provide some additional results deriving from Theorems \ref{relax-1} and \ref{relax0}.

We  recall that for any function $F:X \to \overline{\R}$ the sublevel sets of the relaxed function $\Gamma_{\tau} (F)$ can be represented  by using the following identity that holds for any $\lambda \in \R$
%\begin{equation}\label{sopra}
$$
\{x\in X:\, \Gamma_{\tau} (F)\le \lambda \}=\bigcap_{\lambda'>\lambda}\overline{\{x\in X: \, F(x)\le \lambda'\}}^\tau
$$%\end{equation}
(see Proposition 3.5 in \cite{DM93}). 
Under the hypoyhesis that $\Gamma_{\tau} (F)$ is level convex,  we provide a refined relationship between sublevel sets of $\Gamma_{\tau} (F) $ and $F$  together with 
an abstract representation formula for $\Gamma_{\tau} (F)$.

\begin{prop}\label{reprsublev}  Assume that $\Gamma_{\tau} (F)$ is level convex, then for any $x\in X$ it holds 
\begin{equation}\label{formula}
\Gamma_{\tau} (F)(x)=\inf \left\{\lambda \, :\,  x\in \tau\hbox{-}\conv(\{x'\in X:\,  F(x')\le \lambda\})\right\}, 
\end{equation}
where for $A\subseteq X$ $\tau\hbox{-}\conv(A)$ is defined in \eqref{tauco}.  \\
\noindent In addition, for any $\lambda \in \R$ it holds 
\begin{equation}\label{formula2}
\{x\in X:\,  \Gamma_{\tau} (F)(x)\le \lambda \}=\bigcap_{\lambda'>\lambda}\tau\hbox{-}\conv(\{x\in X:\,  F(x)\le \lambda'\}). 
\end{equation} 
  \end{prop}
 \proof 
Set 
$$
G(x)=\inf \{\lambda: x\in \tau\hbox{-}\conv( \{x'\in X:\,  F(x')\le \lambda\}) \}, $$ 
it can be easily verified that 
$$
\{x\in X:\,  G(x)\le \lambda\}= \bigcap_{\lambda' > \lambda} \tau\hbox{-}\conv (\{x\in X:\,  F(x)\le \lambda'\}).
$$
Hence the thesis follows once we prove that $\Gamma_{\tau} (F)=G$. 
\noindent  To this aim we note the following facts:
\begin{enumerate}
\item[(1)] for every function $f:X\to \bar \R$ it holds 
 $f(x)=\min\{\lambda: x\in \{f(x')\le \lambda\}\};$
 \item[(2)]  for any family of subsets of $X$,  $(C_\lambda)_\lambda$, with $C_\lambda\subseteq C_{\lambda'}$ for  $\lambda\le \lambda'$ and $\bigcup_\lambda C_\lambda=X$ one may define  
 $$f_C(x):=\inf\{\lambda: x\in C_\lambda\}.$$  
%Then  $$
%\{x\in X:\,  f_C(x)\le \lambda\}= \bigcap_{\lambda' > \lambda}  (\{x\in X:\,  x\in C_{ \lambda'}\}).
%$$
If  $(D_\lambda)_\lambda$ is another family of subsets of $X$,  with the same properties of $(C_\lambda)_\lambda$, then 
$$
C_\lambda\subseteq D_\lambda \ \forall\lambda \in \R \Rightarrow f_D(x):=\inf\{\lambda: x\in D_\lambda\}\le f_C(x)\hbox{ for any }x\in X.
$$
\end{enumerate}

%\noindent To this aim we note that for a given function $f:X\to \bar \R$ it holds 
% $$f(x)=\color{red}{\min}\{\lambda: x\in \{f\le \lambda\}\}$$ 
%and for any family of subsets of $X$,  $(C_\lambda)_\lambda$, with $C_\lambda\subseteq C_{\lambda'}$ for  $\lambda\le \lambda'$ and $\cup_\lambda C_\lambda=X$ one may define  
% $$f_C(x)=\inf\{\lambda: x\in C_\lambda\}.$$  
%Moreover,  given  $(D_\lambda)_\lambda$ another family of subsets of $X$,  with the same properties of $(C_\lambda)_\lambda$ we have  
%$$
%C_\lambda\subseteq D_\lambda \ \forall\lambda \in \R \Rightarrow f_D(x):=\inf\{\lambda: x\in D_\lambda\}\le f_C(x)\hbox{ for any }x\in X.
%$$

In order to apply the arguments above to $G, F$ and $\Gamma_{\tau}$, we note that 
 $\tau\hbox{-}\conv (\{F\le \lambda\})\supseteq \{F\le \lambda\}$ for any $\lambda \in \R$. Thus,  by (1) e (2) above,  we deduce that  $G(x)\le F(x)$ for every $x\in X$. Moreover,  by the closedness of $\tau\hbox{-}\conv (\{x\in X:\, F(x)\le \lambda'\})$ for any $\lambda'\in \R$ we get that $G$ is $\tau$-lower semicontinuous. Hence $G\le \Gamma_{\tau} (F)$. On the other hand, since $\{x\in X:\, \Gamma_{\tau} (F)\le \lambda \}$ is a closed convex set by hypothesis and 
$\{x\in X:\,  \Gamma_{\tau} (F)(x)\le \lambda \}\supseteq \{x\in X:\,F(x)\le \lambda\}$, we deduce  
$$
\{x\in X:\,  \Gamma_{\tau} (F)(x)\le \lambda \}\supseteq \tau\hbox{-}\conv\{x\in X:\, F(x)\le \lambda\}.$$
 Hence,  we get  $\Gamma_{\tau} (F)\le G$  and this yields \eqref{formula} and \eqref{formula2}. 
 
 \qed
 
 In the supremal case, once the level convexity of $\Gamma_{\tau}(F)$ is established in Theorems \ref{relax-1} or  \ref{relax0},  one can easily prove the following result. 

\begin{cor}\label{rappr2}  Under the same hypotheses of Theorem \ref{relax-1} or Theorem \ref{relax0},  for any $u\in \wi$ it holds 
\begin{equation}\label{formulasup}
\Gamma_{\tau} (F)(u)=\inf \left\{\lambda \, :\,  u\in \tau\hbox{-}\conv(\{v\in \wi \,:\,F(v)\le \lambda\})\right\} 
\end{equation} 
 and, for any $\lambda \in \R$, we have 
\begin{equation}\label{formula2sup}
\{v\in \wi: \,\Gamma_{\tau} (F)(v)\le \lambda \}=\bigcap_{\lambda'>\lambda}\tau\hbox{-}\conv(\{v\in \wi \,:\,F(v)\le \lambda'\}),  
\end{equation} 
when $\tau$ is one of  the topologies $\tau_{\infty}, w^*, w^*_{seq}$. 
%\st{ In addition, it holds $\Gamma_{w^*} (F)=\Gamma_{w^*_{seq}} (F)$.} {\color{red} \bf Toglierei anche io resto. }
%
%$$F \hbox{ is w*-lower semicontinuous } \Longleftrightarrow F  \hbox{ is sequentially w*-lower semicontinuous.}$$

\end{cor} 

\proof By the level convexity of  $\Gamma_{\tau}(F)$, by applying Proposition~\ref{reprsublev}, we deduce \eqref{formulasup} and \eqref{formula2sup}. %Finally the last part of the statement follows by Proposition \ref{weakseq}. 
 
 \qed
 
\begin{rem} {\rm We note that, under hypothesis \eqref{dFbzero}, taking into account \eqref{conlevelset} in Proposition \ref{chiusura}, we can also obtain that 
$$
\{v\in \wi: \,\Gamma_{\tau} (F)(v)\le \lambda \}=\bigcap_{\lambda'>\lambda}\tau\hbox{-}\conv\Big(\{ u\in\wi : \  F(u)\leq\lambda  \}\Big).
$$}
\end{rem}

We conclude this section  with a representation result which improves \cite[Theorem 2.3]{Pr08} since we do not require a continuity assumption on $f(x,\cdot)$. Note that, as shown in Remark 3.1 of [14], the level convexity of the supremal functional $F$ does not imply the level convexity of $f(x,\cdot)$. However, it is possible to show  that any level convex supremal functional F can be represented by a level convex supremand $\f$ (possibly different from the level convex envelope of $f(x,\cdot)$, see \cite[Example 8.1]{Pr08}).

In order to show the result, we recall that  a Borel function  $f: \Omega \times \R^N \to \bar \R$ is a  {\sl normal supremand} if  for a.e. $x\in\Omega$ the function $f(x,\cdot)$ is lower semicontinuous on $\R^{ N}$.

\begin{prop}\label{equivalenza} Let  $\Omega$ be a bounded open set   with Lipschitz continuous boundary. Let  $f:\Omega\times \R^N\to \bar\R$ be a normal supremand.  Let $F$ be  the    supremal functional (\ref{sfzero}) represented by  $f$ and  assume that $F$ satisfies the  condition  \eqref{accaomega}. Then the following facts are equivalent:
\begin{enumerate}
\item [(i)] $F$ is  $w^*$-lower semicontinuous in   $\wi$;
\item [(ii)]$F$ is  $w_{seq}^*$-lower semicontinuous  in  $\wi$;
\item [(iii)]  $F$ is a level convex supremal functional;
 \item [(iv)]  there exists a  level convex normal supremand $\f:\Omega\times \R^N\to \overline{\R}$ given by
\begin{equation}\label{formulavarphi}\f(x,\xi):=\inf\left\{  F(u)\; |\; \begin{array}{l} u\in W^{1,\infty}(\om)\;
\hbox{ s.t.  $x\in \widehat{u},$ with $Du(x)=\xi$}\end{array}\right\}
\end{equation}
where $$\widehat{u}:=\left\{x\in\Omega : \mbox{\rm $x$ is a Lebesgue point of $Du$ and a differentiable point of $u$} \right\}$$
such that  $$
F(u) = \supess_{x\in \Om} \f(x, D u(x)).
$$
Moreover  there esists a negligible set $N\subset \Omega$ such that $\f(x,\xi)\geq f(x,\xi)$  $\forall\, x\in\Om\setminus N$, $\forall\, \xi\in \R^N$.

\end{enumerate}
\end{prop}

\proof (i) $\Longrightarrow$ (ii) follows by Remark \ref{lsc}(2) while (ii) $\Longrightarrow$ (iii) follows by Theorem \ref{relax0}. 
In order to prove  (iii) $\Longrightarrow$ (iv) let us consider  the auxiliary functional $G:\wi\to [0,+\infty]$ given by $$G(u):= \arctan F(u)+\frac \pi 2=\supess_{\Omega} \arctan  (f(x,Du(x))+\frac \pi 2.$$ Note that  $G$ is level convex if and only if $F$ is level convex.
Thanks to the result  \cite[Theorem 2.4]{Pr08}, the level convex functional $G$ can be represented through the level convex normal supremand $g$ given by 
$$g(x,\xi):=\inf\left\{ G(u)\; |\; \begin{array}{l} u\in W^{1,\infty}(\om)\;
\hbox{ s.t.  $x\in \widehat{u},$ with $Du(x)=\xi$}\end{array}\right\}
$$ and satisfying  $$g(x,\xi)\geq \arctan f(x,\xi)+ \frac \pi 2 \quad \forall\, x\in\Om\setminus N\,, \forall\, \xi\in \R^N$$
where $N\subset \Omega$ is a  negligible set.
Hence it is sufficient to choose $$\f(x,\xi):=\tan (g(x,\xi)- \frac \pi 2)= \inf\left\{ F(u)\; |\; \begin{array}{l} u\in W^{1,\infty}(\om)\;
\hbox{ s.t.  $x\in \widehat{u},$ with $Du(x)=\xi$}\end{array}\right\}. $$ 
Eventually, (iv) $\Longrightarrow$ (i) follows by \cite[Theorem 3.4]{BJW99}.

\qed

\section{Some interesting examples and counterexamples}\label{examples} 

In this section we collect some examples pointing out the optimality of some of  the statements given in the previous sections. 

\subsection{An infimum of supremal functionals that is not supremal}

Given  $F,G:\wi\to \bar\R$  supremal functionals represented by the supremands $f,g:\Omega \times \R^N\to \bar\R$, respectively we have seen in Example \ref{specialcase} that the functional  $F\vee G$  defined by $F\vee G(u):=F(u)\vee G(u)$ is still a supremal functional. On the other hand, this property may not hold in general for  the functional $F\wedge G(u):=F(u)\wedge G(u)$, as shown in the following example. 
\begin{ex}{\rm  Let $\Omega=(-1,1)$ and let $f,g:(-1,1)\to \R$ be defined as 
$$f(x):=\begin{cases}
1 &  \hbox{if } x\in (-1,0)\cr
\cr 
3 &  \hbox{if } x\in (0,1)\cr

 \end{cases}
 , 
\quad g(x):=\begin{cases}
4 &  \hbox{if } x\in (-1,0)\cr
\cr
2 &  \hbox{if } x\in (0,1)\cr 
 \end{cases}
.
$$

\noindent We  consider the (localized) supremal functionals $F,G$ with supremands $f,g$, respectively, given by  $F(u,A)=\supess_{A}f(x)$ and $G(u,A)=\supess_{A}g(x)$ for any open set $A\subseteq (-1,1)$ and $u\in W^{1,\infty}((-1,1))$. We claim that $F\wedge G$ cannot be represented in a supremal form since it does  not satisfy the necessary condition $$
F\wedge G(u,\bigcup_{i\in I} A_i)=\bigvee_{i\in I} (F\wedge G)(u,A_i) \quad \forall\, u\in \wi,\ \forall A_i\in {\mathcal A}(\Omega).
$$
Indeed, set $A=(-1,-\frac 12)$ and $B=(0,1)$, an easy computation shows that 
$$
F\wedge G (u, A\cup B)= 3> 2=(F\wedge G(u, A))\vee (F\wedge G(u,B)).
$$
}
\end{ex}

\subsection{A supremal functional $F$ with  $d^\lambda_F <d^\lambda_{\Gamma_\tau (F)}$ for some $\lambda$}

In the following example we show that, given  a supremal functional $F$, the intrinsic distance associated to its relaxed functional $\Gamma_\tau (F)$ is in general different from  $d^\lambda_F$, for some values $\lambda >\inf_{\wi} F$.

\begin{ex}\label{gap} {\rm 
Let $\Omega=(0,1)$ and $F(u)=\supess_{\Omega} f( u' (x))$ where $f:\R\to \R$ is given by 
$$f(z):=\begin{cases}
-z &  \hbox{if } z\le 0\cr
\cr 
1+\dfrac{1}{z} &  \hbox{if } z>0.\cr

 \end{cases}
$$
It can be verified that, for any $\tau\in \{ \tau_{\infty}, w^*, w^*_{seq}\}$,  the relaxed functional $\Gamma_\tau (F)$ is still a supremal functional given by 
$\Gamma_\tau (F)(u)=\supess_{\Omega} \phi( u' (x))$ where $\phi:\R\to \R$ is given by 
$$\phi(z):=\begin{cases}
-z &  \hbox{if } z\le 0\cr
\cr 
1 &  \hbox{if } z>0.\cr

 \end{cases}
$$
Taking aside the computation above we will deduce the value of the intrinsic distance functions by computing the level sets of $F$  and taking advantage of the results established in Corollary \ref{rappr2}.   Indeed, the sequence  $u_n$ given by $u_n(x)=-\frac{x}{n}$ is a minimizing sequence along which $F$ is continuous and thus $F$ satisfies \eqref{accaomega}. 

A direct computation shows that 
\begin{equation}\label{livello1}
\{u\in \wi\,:\, F(u)\le 1\}= \{u\in \wi\,:\, -1\le u'(x)\le 0\hbox{ a.e. in } \Omega\}
\end{equation}
and, for a given $\lambda>1$, we have also 
$$
\{u\in \wi\,:\, F(u)\le \lambda\}= \{u\in \wi\,:\,  u'(x)\in [-\lambda, 0]\cup [1/(\lambda-1), +\infty) \hbox{ a.e. in } \Omega\}.
$$
By \eqref{livello1} we deduce that   
$$d^1_F(x,y)=\begin{cases}
0 &  \hbox{if } x\ge y\cr
\cr 
|x-y| &  \hbox{if } x<y.\cr

 \end{cases}
$$
On the other hand, it can be shown that for any $\tau\in \{ \tau_{\infty}, w^*, w^*_{seq}\}$ it holds 
 $$
\tau\hbox{-}\conv (\{u\in \wi\,:\, F(u)\le \lambda\})= \{u\in \wi\,:\,  u'(x)\in [-\lambda,  +\infty) \hbox{ a.e. in } \Omega\}.
$$
By applying formula \eqref{formula2sup} in Corollary \ref{rappr2} we deduce that 
$$
\{u\in \wi\,:\, \Gamma_\tau(F)(u)\le 1\}= \{u\in \wi\,:\, -1\le u'(x)\hbox{ a.e. in } \Omega\}
$$
and 
$$d^1_{\Gamma_\tau (F)}(x,y)=\begin{cases}
+\infty &  \hbox{if } x\ge y\cr
\cr 
|x-y| &  \hbox{if } x<y.\cr

 \end{cases}
$$
Note that this example also provides a counterexample to the validity of the formula 
$$
\{v\in \wi: \,\Gamma_{\tau} (F)(v)\le \lambda \}=\bigcap_{\lambda'\ge\lambda}\tau\hbox{-}\conv(\{v\in \wi \,:\,F(v)\le \lambda'\}).
$$
 
}\end{ex}

\subsection{Functional $F$ with discontinuous supremand satisfying \eqref{accaomega}.} 
In the following example we exhibit a supremal functional $F$ represented by a discontinuos supremand $f$ 
whose minimum is not attained. We construct a "minimizing" sequence of $F$ made up by  "continuity" points  which satisfies \eqref{accaomega}.

\begin{ex} {\rm  Let $F:W^{1,\infty}(0,1)\to \R$ be a functional whose supremand is given by 
 
$$f(\xi):=\begin{cases}
\xi +2  &  \hbox{if } \xi\ge 0\cr
\cr 
-\xi &  \hbox{if } \xi<0.\cr

 \end{cases}
$$
For such functional \eqref{dFbzero} fails as it can be easily shown that  $d^1_F(x,y)=0$ if $x<y$. Despite of this  if we set $u_n(x):= -\dfrac{x}{n+1}$, $n\in \N$, then  it holds 
$$ \lim_{\xi\to 0} F(u_n+u_\xi)=F(u_n), \quad 
\lim\limits_{n\to +\infty}F(u_n)=0=\inf_{\wi} F.
$$
In particular \eqref{accaomega} is satisfied.

}\end{ex}

Here below we give an example of a supremal functional $F$  represented by a function  $f$ that is discontinuous everywhere in $\R$ and still satisfying \eqref{accaomega}.

\begin{ex} {\rm 
 We consider the following function $:[0,1]\to \R$ 
 \begin{equation}g(\xi):=\begin{cases} 0 &  \hbox{if } \xi\in [\frac 1 {n}, \frac 1 {n-1})\cap \Q,n\in \N,n\geq 2\cr
 \\
 
\frac 1 {n}&  \hbox{if } \xi\in [\frac 1 {n}, \frac 1 {n-1} )\; \cap (\R\setminus \Q),\, n\geq 2\cr
\\
1  & \hbox{ otherwise}
. \end{cases} 
\end{equation}
 Moroever for every $ \la \in [\frac 1{ n}, \frac 1 {n-1})$  we have that

\begin{equation}\{ g(\xi)\leq \la\} =(0, \frac 1 {n-1}]\cup \left( \Q \cap (\frac 1 {n-1},1) \right).
\end{equation}
The function $g$ is discontinuous everywhere in $[0,1]$. In particular, its periodic extension  $f$ on $\R$ is discontinuous everywhere in $\R.$
 However  the supremal functional  $F(u)=\supess_{x\in \Omega} f(Du(x))$ satisfies \eqref{accaomega} since the sublevel set  of  $f$ have not empty interiors.
% Note that there don't exist $(u_n)_{n\in \N}\subseteq \wi$ and a sequence $(\la_n)\to 0$ such that 
%$$\begin{cases}
%F(u_n)\to 0\\
% \lim\limits_{\xi\to 0} F(u_n+u_\xi)=F(u_n) \quad  \forall\, n\in \N
% \\
% \end{cases}
%$$
%In fact if $F(u_n)\to 0$ then $F(u_n)\leq \frac 1 k$ for $n$ big enough.\\
% Therefore $u'_n(x)\in ]0, \frac 1 {k-1}]\cup ((\frac 1 {k-1},1[\cap \Q)$ for a.e.$x$.\\
%
%If $u'_n(x)\in \Q$ for a.e. $x$ then $f(u'_n(x)+\xi)=\frac 1 k$ for all $\xi\in \R\setminus \Q$ small enough while
%$f(u'_n(x)+\xi)=0$ for all $\xi\in \Q$ small enough.
%Then $F(u_n+u_\xi)\not \to 0$ when $\xi\to 0$.\\
%
%if $u'_n(x)\in \R\setminus \Q$ for a.e. $x$ then $f(u'_n(x)+\xi)=\frac 1 k$ for all $\xi\in  \Q$ small enough while
%$f(u'_n(x)+\xi)=0$ for all $\xi\in \Q$ small enough.
%Then $F(u_n+u_\xi)\not \to 0$ when $\xi\to 0$.
 }\end{ex}

\subsection{Failure of the local representation by means of the supremand $\f$}
In the following example we show that, given a normal supremand  $f$, the function $\f$ in  the representation formula \eqref{formulavarphi} is in general affected by the reference set $\Omega$. More in details  we exhibit a supremal functional $F$ such that the supremand $\f$ does not represent the lower semicontinuous envelope of the localized version of $F$,  $F(\cdot,A)$ defined by  
\begin{equation}\label{localiz}
F(u,A) := \supess_{A} f(x, D u(x)),
\end{equation} 
for $A\in \mathcal{A}(\Omega)$.

\begin{ex}\label{boh}{\rm  Let $\Omega=(-2,2)$ and let $f:(-2,2)\times \R\to \R$ be defined as 
$$f(x, \xi):=\begin{cases}
(1 -|\xi|)\vee 0 &  \hbox{if } x\in [-1,1]\cr
\cr 
2+|\xi| &  \hbox{if } x\in (-2,-1)\cup(1,2) .\cr

 \end{cases} 
$$
\noindent Let $F$ be the supremal functional given by  $F(u)=\supess_{\Omega}f(x, u'(x))$.  Note that $F$ is level convex on $\wi$ since $F(u)= 2+\supess_{(-2,-1)\cap (1,2)} |u'(x)|$ for every $u\in \wi$. 

A direct  computation shows that the function $\varphi$ in Corollary \ref{equivalenza}(iv) is given by 
$$
\f (x, \xi)=\begin{cases}
2 &  \hbox{if } x\in (-1,1)\cr
\cr 
2+|\xi| &  \hbox{if } x\in (-2,-1)\cup(1,2). \cr
\end{cases}
$$
\noindent  Note that $\f$  is level convex and is strictly greater than $f$ on a set of positive measure.  

\noindent We underline that the function $\f$ is obtained by a minimization process on the whole set $\Omega$ and it is not suitable to represent the lower semicontinuous envelope of the localized functional $F(\cdot,A)$ as shown by the following computation. For any open set $A\subseteq (-2,2)$ let  $F(u,A):=\supess_{A}f(x, u'(x))$ be the supremal functional defined on  $W^{1,\infty}((-2,2))$ and let $G(u,A):=\Gamma_\tau (F) (u, A)$ where  $\tau$ is one of the  topologies $\tau_{\infty}, w^*, w^*_{seq}$. We have that $f$ is a normal supremand and  $f(x,\cdot)$ is level convex in $\R$ if and only if  $x\in \Omega\setminus [-1,1]$. According to this the functional $F(\cdot, A)$ is level convex  only on the open sets $A$ satisfying the condition $|A\setminus (-1,1)|>0$.  Indeed, for such an open set  $A$ it holds 
$$
F(u,A) = \supess_{A\setminus (-1,1)} f(x,  u'(x))
$$
for any $u\in \wi$ and the supremand $f(x, \cdot)$ is level convex for any $x\in \Omega\setminus [-1,1]$.

\noindent We claim that for any open set $A$  the relaxed functional $G (\cdot, A)$ can be represented  by the function $g:(-2,2)\times \R\to \R$ defined as 
$$g(x, \xi):=\begin{cases}
0 &  \hbox{if } x\in (-1,1)\cr
\cr 
2+|\xi| &  \hbox{if } x\in (-2,-1)\cup(1,2). \cr

 \end{cases}
$$
Indeed, if $A$ is an open set in $\Om$ with $|A\setminus (-1,1)|>0$ it holds 
$$
F(u,A) = F(u,A\setminus [-1,1])= \supess_{A\setminus [-1,1]} f(x,  u'(x))
$$
for any $u\in \wi$ and $G(u,A)=F(u,A)$ since $F(u,A\setminus [-1,1])$ is lower semicontinuous in $\wi$ by the level convexity of $f$.
\noindent It remains to prove that for any fixed $A$ with  $A\subseteq (-1,1)$ it holds  $G(u, A)=0$ for any $u\in \wi$. 
Indeed, let $u\in \wi$ be fixed and set $C=\supess_A |u'(x)|$. For $n\in \N$ let us define $\psi_n$  as the odd function piecewise affine such that $\psi_n(0)=0$ and $$\psi'_n(x):= \begin{cases} C+1 & x\in (\frac{2k}{2n},\frac{ 2k+1}{2n})\cr
-(C+1) & x\in  (\frac{2k+1}{2n},\frac{ 2k+2}{2n})\cr \end{cases}$$ with $k=0, \ldots [\frac{n}{2}]$.  We have that the sequence $u_n=u+\phi_n$ converges to $u$ in any topology above and $|u'_n(x)|>1$ for any $x\in A$. Hence $F(u_n, A)=0$ and, subsequently, $\Gamma_\tau(F)(u,A)=0$.   As a consequence, for every $A\in \mathcal{A}(\Omega)$ we have verified that 
$$
G(u,A) = \supess_{A} g(x, u'(x)) .
$$
In particular, for  $A\subset (-1,1)$, we have that $G(u,A) <\supess_{A} \f(x, u'(x))$ for any $u\in\wi$. }
\end{ex}


\begin{thebibliography}{CTCR81}




\bibitem{AP} { N. Ansini, F. Prinari}. {\it Power law approximation of supremal functional under differential constraint.}  SIAM J. Math. Anal. {\bf 2} 46 (2014),  
1085--111.
%\bibitem{AP2} {\sc N. Ansini, F. Prinari} {\it ''Lower semicontinuity of supremal functional under differential constraint''}  { ESAIM Control Optim. Calc. Var.}  21 (4) (2015), pp. 1053-1075.


\bibitem{BJW99}{ E.N. Barron, R. R. Jensen, C.Y.Wang}. {\sl
Lower Semicontinuity of \( L^{\infty }\)-Functionals.} Ann. Inst. H.
Poincar\'e Anal. Non Lin\'eaire (4) {\bf18} (2001), 495--517.

\bibitem{BJW}{ E.N. Barron, R. R. Jensen, C.Y.Wang}. {\sl
 The Euler Equation and Absolute Minimizers of \( L^{\infty }\)-Functionals}. Arch. Rational Mech. Anal.  157 (2001), 225--283.

\bibitem{BN}{  M. Bocea and V. Nesi}. {\sl $\Gamma$-convergence of power-law functionals, variational principles in  $L^\infty$ and applications}. SIAM J. Math. Anal. {\bf 39} (2008), 1550--1576.

\bibitem{Brezis}{ H. Brezis}: {\sl Functional Analysis, Sobolev Spaces and Partial Differential Equations.}  Universitext, Springer, New York (2011).

\bibitem{BGP}{ A. Briani, F. Prinari, A. Garroni}. {\sl  Homogenization of $L\sp \infty$-functionals.}   Math. Models Methods Appl. Sci.  14 (2004), no. 12, pp. 1761--1784.
 

\bibitem{Bu89}{ G.  Buttazzo}. {\sl Semicontinuity, Relaxation and Integral Representation in the
 Calculus of Variation.} Pitman Research Notes in Mathematics Series {\bf 207},
Harlow (1989).
 

\bibitem{CDPP02} { T. Champion, L. De Pascale, F. Prinari}. {\sl Semicontinuity
  and absolute minimizers for supremal functionals.} 
   ESAIM Control Optim. Calc. (1) {\bf 10}, (2004), 14--27.

\bibitem{CDP} { T. Champion, L. De Pascale}. {\sl Principles of comparison with distance functions for AML.} J. Convex Anal. 14 (2007), no. 3, 515--541. 

\bibitem{DM93} {G. Dal Maso}. {\sl An Introduction to \( \Gamma\)-Convergence.} Progress in Nonlinear
Differential Equations and their Applications {\bf 8}, Birkhauser, Boston
(1993).

\bibitem{DP} { A. Davini, M. Ponsiglione}.  {\sl Homogenization of two-phase metrics and applications. } 
{J. Anal. Math.}, 103 (2007),  157--196.


 \bibitem{DCP1} {G. De Cecco, G. Palmieri} {\sl Integral distance on a Lipschitz Riemannian manifold.} Math. Z., 207 (1991), 223-243.
\bibitem{DCP2} {G. De Cecco, G. Palmieri} {\sl LIP manifolds:from metric to Finslerian structure.} Math. Z., 218 (1995), 223-237.
 \bibitem {GPP}{ A. Garroni, M. Ponsiglione, F. Prinari}.  {\sl From $1$-homogeneous supremal functionals to difference quotients:
relaxation and $\Gamma$-convergence.} 
{Calc. Var. Part. Diff. Eq}, 27 (2006), no. 4,  397--420.

\bibitem{MG} { M. Gori, F. Maggi.} {\it  On the lower semicontinuity of supremal functionals.} ESAIM Control Optim. Calc. Var. {\bf 9} (2003), 135--143. 

\bibitem{KSZ} { P. Koskela, N. Shanmugalingam, Y. Zhou.} {\it  Intrinsic geometry  and analysis of diffusion processes and $L^\infty$-variational problems.} Arch. Rational Mech. Anal.  214 (2014), 99--142. 
\bibitem{N}{J. R. Norris} {\it Heat kernel asymptotics
and the distance function in Lipschitz Riemannian manifolds} Acta Math., 179 (1997), 79-103

\bibitem{Pr08}{ F. Prinari}.  {\it Semicontinuity and supremal representation in Calculus of Variations}, Appl. Mat. Optim., {\bf 58}, (2008), 111--145.

 
 \bibitem{P09}{ F. Prinari.} {\it  Semicontinuity and relaxation of $L^{\infty}$-functionals.}
  Adv. Calc. Var. (1) {\bf 2} (2009), 43--71.


\end{thebibliography}
\end{document}